\documentclass{article}
\usepackage{graphicx,hyperref}
\usepackage{epsfig,latexsym,xspace,subfigure,color,todonotes}
\usepackage{oupbib}

\newcommand{\citep}{\cite}
\newcommand{\citet}{\textcite}

\newcommand{\addo}[2]{ {\color{blue}#2}} \newcommand{\takeo}[2]{
  {\color{red}\sout{#2}}}
\newcommand{\noteo}[2]{\todo[inline,color=orange!50]{{{#1} note
      start}}{\noindent\color{orange}#2}\todo[inline,color=orange!50]{\hfill{#1}
    note end}\noindent}

\newcommand{\addm}[1]{\addo{Monica}{#1}}
\newcommand{\takem}[1]{\takeo{Monica}{#1}}

\newcommand{\addp}[1]{\addo{Philip}{#1}}
\newcommand{\takep}[1]{\takeo{Philip}{#1}}
\newcommand{\notep}[1]{\noteo{Philip}{#1}}

\newcommand{\adds}[1]{\addo{Steve}{#1}}
\newcommand{\takes}[1]{\takeo{Steve}{#1}}
\newcommand{\notes}[1]{\noteo{Steve}{#1}}

\renewcommand{\addo}[2]{#2} \renewcommand{\takeo}[2]{}
\renewcommand{\noteo}[2]{}

\newcommand{\p}{\mbox{P}\xspace} \newcommand{\pc}{\mbox{PC}\xspace}
\newcommand{\alte}{\mbox{ALTE}\xspace}
\newcommand{\bleed}{\mbox{bleed}\xspace}
\newcommand{\abuse}{\mbox{abuse}\xspace}
\newcommand{\rr}{\mbox{RR}\xspace}
\newcommand{\orr}{\mbox{ORR}\xspace}

\newcommand{\ie}{{\em i.e.\/}\xspace}
 \newcommand{\cf}{{\em
    cf.\/}\xspace} 
\newcommand{\E}{{\mbox{E}}}

\newcommand{\secref}[1]{\mbox{\S$\,$\ref{sec:#1}}}
\newcommand{\eqref}[1]{\mbox{(\ref{eq:#1})}}
\newcommand{\tabref}[1]{\mbox{Table~\ref{tab:#1}}}
\newcommand{\figref}[1]{\mbox{Figure~\ref{fig:#1}}}
\newcommand{\winbugs}{{\sc WinBUGS}\xspace} \newcommand{\eg}{{\em
    e.g.\/}\xspace} \newcommand{\cd}{\,|\,}
\newcommand{\cip}{\mbox{$\perp\!\!\!\perp$}}

\newtheorem{cond}{Condition}[section]
\newcommand{\condref}[1]{\mbox{Condition~\ref{cond:#1}}}

\renewcommand{\thefootnote}{\fnsymbol{footnote}}

\author{A.~P.~Dawid\thanks{University of Cambridge} \and
  M.~Musio\thanks{University of Cagliari} \and
  S.~E.~Fienberg\thanks{Carnegie Mellon University}}
\title{From Statistical Evidence to Evidence of Causality}

\begin{document}
\maketitle 

% \footnotetext[1]{University of Cambridge
%   \href{mailto:apd@statslab.cam.ac.uk}{apd@statslab.cam.ac.uk} }
% \footnotetext[2]{University of Cagliari
%   \href{mailto:mmusio@unica.it}{mmusio@unica.it}}
% \footnotetext[3]{Carnegie Mellon University
%   \href{mailto:fienberg@stat.cmu.edu}{fienberg@stat.cmu.edu} }

\renewcommand{\thefootnote}{\arabic{footnote}}

\nocite{%national2011Reference,
  studies}

\begin{abstract}
  \takem{Science is largely concerned with understanding the ``effects
    of causes'' (EoC), while Law is more concerned with understanding
    the ``causes of effects'' (CoE)}\addm{\addp{\noindent While
    }statisticians and quantitative social scientists \takep{are most
      comfortable with thinking about}\addp{typically study} the
    ``effects of causes'' (EoC), \takep{while }Lawyers \adds{and the
      Courts} are more concerned with understanding the ``causes of
    effects'' (CoE)}.  \takep{While }EoC can be addressed using
  experimental design and statistical analysis, \addp{but }it is less
  clear how to incorporate statistical or epidemiological evidence
  into CoE reasoning, as might be required for a case at Law.  Some
  form of counterfactual reasoning, such as \takep{Rubin's ``potential
    outcomes'' approach}\addp{the ``potential outcomes'' approach
    championed by Rubin}, appears unavoidable, but this typically
  yields ``answers'' that are sensitive to arbitrary and untestable
  assumptions.  We must therefore recognise that a CoE question simply
  might not have a well-determined answer.  It \takep{may}\addp{is}
  nevertheless \takep{be }possible to use statistical data to set
  bounds within \addm{which} any answer must lie.  With less than
  perfect data these bounds will themselves be uncertain, leading to a
  \takem{novel} compounding of \takep{two}\addp{different} kinds of
  uncertainty.  Still further care is required in the presence of
  possible confounding factors.  In addition, even identifying the
  relevant ``counterfactual contrast'' may be a matter of Policy as
  much as of Science.  Defining the question is as non-trivial a task
  as finding a route towards an answer.  This paper develops some
  technical elaborations of these philosophical points, and
  illustrates them with an analysis of a
  case study in child protection.\\

\end{abstract}

\section{Introduction}
\label{sec:intro}

One function of a Court of Law is to attempt to assign responsibility
or blame for some undesirable outcome.  In many such cases there will
be relevant testimony about statistical or epidemiological evidence
arising from studies done on specialized populations, but this
evidence addresses the main issue only indirectly, at best.  It has
until now been unclear how to use such evidence to focus on the issue
at hand which involves specific individuals experiencing the
undesirable outcome.  Although there is a considerable literature on
certain aspects of this problem---see for example \citet{refepi},
which aims to assist US judges in managing cases involving complex
scientific and technical evidence---we consider that there are
important logical subtleties that have not as yet been accorded the
appreciation they warrant.  Here we show that, even in the (very rare)
case that we have the best possible and most extensive data \adds{on
  the ``effects of causes''}, and can accept certain very strong but
necessary conditions, there will still remain irreducible
uncertainty---which we can express as interval bounds---about the
relevant ``probability of causation''.  With less than fully perfect
data, this interval uncertainty will be further compounded by
statistical uncertainty. Such multiple \adds{forms of} uncertainty
\takep{raises new questions}\addp{raise subtle issues} of
interpretation and presentation.

The structure of the paper is as follows.  In \secref{frachon} we
consider a high-profile case where serious side-effects led to the
withdrawal of a drug from the market, and, in turn, to litigation
against the manufacturer.  \adds{Since the evidence in this litigation
  has yet to be presented formally at trial in court, we} \takes{We}
consider how general evidence of incidence of effects might or might
not be relevant to a hypothetical tort action in which an affected
patient sues the manufacturer for damages, and we relate this to the
distinction we draw in \secref{coeeoc} between ``effects of causes''
and ``causes of effects''.  After a brief consideration of inference
from statistical data about effects of causes in \secref{stateoc}, the
remainder of the paper focuses on inference about causes of effects,
based on a ``probability of causation'' defined using counterfactual
logic.  Although this probability is typically impossible to pinpoint
on the basis of epidemiological data, however extensive, \addp{in
  \secref{coe}} we give bounds between which it must lie---bounds
which, however, will themselves be subject to statistical
uncertainty\addp{, which we discuss in \secref{statunc}}.  In
\secref{best} and \secref{datanal} we illustrate our theory with a new
analysis of a case study in child protection.
\addp{Section~\ref{sec:conc} offers some concluding remarks.}

\section{Epidemiological Evidence in Litigation}
\label{sec:frachon}

\subsection {Epidemiological background}
The drug Mediator, also known as benfluorex, was for many years
marketed as an anti-diabetic drug.  It was also widely used off-label
as an appetite suppressant.  In November 2010, however, following the
publication of a popular book by Ir\`ene \citet{frachon:book}, the
French Health Agency CNAM
% (Caisse Nationale d'Ass\-ur\-ance Maladie)
announced its finding that around 500 deaths in France over a thirty
year period could be attributed to Mediator---see also
\citet{Hill:2011}.  This was based on extrapolation of results in two
scientific studies, published at about the same time, focusing on the
effects of benfluorex on valvular heart disease.
\citet{frachon:plos1}\ showed a significantly higher prevalence of
unexplained valvular heart disease in patients taking benfluorex, as
compared to controls.  \citet{Weill:2010}\ examined the records of
over a million diabetic patients in a cohort study, and reported a
higher hospitalisation rate for valvular heart disease in benfluorex
takers.

\subsection{Litigation}
As the news about Mediator reverberated through the media, the French
authorities withdrew the drug from sale.  At the same time, hundreds
of individuals jointly filed a criminal lawsuit against the
manufacturer of Mediator, the French pharmaceutical giant \adds{Les
  Laboratoires Servier}\takes{Servier}. The trial has been under way
since May 2012, with initial aspects focused on whether the company
was guilty of misconduct.  \takep{As of}\addp{At} the time of
preparation of this article, the issue of whether Mediator was in fact
the cause of the heart disease in any of those who brought the lawsuit
had yet to be addressed, and no expert scientific testimony had been
presented to the court. \adds{As of September, 2014, however, the
  company has agreed to compensate over 350 individual
  plaintiffs.\footnote{http://www.lejdd.fr/Societe/Sante/Mediator-350-victimes-indemnisees-685286\#new-reactions}}

In the US benfluorex was removed from the marketplace in the
1990s. \takep{But the}\addp{The} banning in 1997 of a related drug,
Redux, led to a \$12 billion settlement, following a class action by
thousands of individuals \citep{mnt}.  Thus considerable attention
both in France and elsewhere is focused on the case against Servier.

\subsection{Scientific results}

The matched case-control study of \citet{frachon:plos1}\ involved 27
cases of valvular heart disease and 54 controls. Investigators
determined whether the patients had or had not used benfluorex.

We display the core data in \tabref{twobytwo}.  The face-value odds
ratio in this table is $(19\times51)/(3\times8) = 40.1$, but this
could be misleading because of confounding factors\addp{\footnote{See
    \citet{HollandRubin1988} for when it is necessary and appropriate
    to adjust for covariate information in such a study.}}.  A
logistic regression analysis reported by \citet{frachon:plos1}
adjusted for body mass index, diabetes and dexfenfluramine use, and
reduced the odds ratio to $17.1$ \addm{($95\%$ CI $3.5$ to $83.0$)}, a
value which is still a large and highly significant measure of
positive association between benfluorex and valvular heart disease.
In the same direction, \citet{Weill:2010} computed a risk ratio
(though with relatively crude adjustments) of \takem{the order of
  3}\addm{$3.1$ ($95\%$ CI $2.4$ to $4.0$).}
\begin{table}[htdp]
  \begin{center}
    \begin{tabular}{|c|cc|c|}\hline
      Benfluorex Use & Cases & Controls & Totals\\ \hline
      Yes & 19 & 3 &22\\
      No & 8& 51 & 59\\ \hline
      Totals & 27 & 54 & 81\\ \hline
    \end{tabular}
  \end{center}
  \caption{Raw results from case-control study linking benfluorex
    and valvular heart disease.   Source:
    {\protect\citet{frachon:plos1}}.}
  \label{tab:twobytwo}
\end{table}

\subsubsection{Robustness of odds ratio}
\label{sec:robust}

While the risk ratio may be a more relevant and incisive measure of
the strength of an effect than the odds ratio (and it will feature
importantly in our analysis of \secref{stateoc} below), it faces a
very serious problem: it is simply not possible to compute it from a
retrospective study, such as that of \citet{frachon:plos1}.  In
contrast, the odds ratio, whether raw or adjusted via a logistic
regression, has the important property that it is simultaneously a
meaningful measure of association \takep{\citep{BFH:1975}} and
computable from retrospective data
\takep{\citep{Farewell:1979}}\addp{\citep{altham1970,BFH:1975,Farewell:1979}}.
Furthermore, \addp{under suitable adjustment for covariates} it will
be a good approximation to the risk ratio when the outcome is rare.

\subsection{Toxic tort---A hypothetical case}
\label{sec:hypo} Now consider a (currently purely hypothetical) case
that might be brought on the basis of these scientific reports.  A
woman with unexplained valvular heart disease sues the manufacturer of
benfluorex, claiming that it was \takep{that}\addp{this} that caused
her illness.  An epidemiologist, testifying \addp{as expert witness}
for plaintiff, claims that, on the evidence of Dr.\ Frachon's and Dr.\
Weill's studies, the medication can cause valvular heart
disease. \adds{\takep{especially since }\addp{Citing
    \citet{Nicot:2011}, who} argued that ``the probabilistic
  information, derived from the available epidemiological studies,
  needs to be considered as part of evidence to establish or refute a
  causal link between benfluorex and valvular disease for a given
  patient''\addp{, this witness goes on to claim that this is evidence
    for a causal link in the current case}}. The defendants in turn
proffer their expert, who testifies that in the manufacturer's
clinical trials there was no evidence of such a side effect.  How
should the court rule?

The court needs to decide on the cause of this woman's heart
disease. But the plaintiff's expert addresses something different, the
general scientific question ``Can benfluorex be shown to cause heart
disease?'' For an epidemiologist, the evidence for this would ideally
be captured by the risk ratio, though, as we have seen, for the
Frachon data we would have to be satisfied with the adjusted odds
ratio instead.  But even if we had perfect and unassailable
statistical evidence in support of this general scientific hypothesis,
that would still be only very indirectly relevant to the individual
case at issue.  We shall see below that the relationship between such
a generalisation and the specific issue before the court is extremely
subtle. \adds{For an extended discussion of the legal nexus of
  individual causation and the ``causes of effects'', see
  \citet{fienberg2014}}.

\section{Causes of Effects and Effects of Causes}
\label{sec:coeeoc}

One might be tempted to assume that the ``effects of causes''
(henceforth EoC) and the ``causes of effects'' (CoE) are related
probabilistically via Bayes theorem.  After all, this was how
\citet{laplace1986}\ introduced the topic: ``If an event can be
produced by a number $n$ of different causes, the probabilities of
these causes given the event are to each other as the probabilities of
the event given the causes,\ldots''\@ Later authors recognized the
issue to be more complex.

John Stuart Mill distinguished between inferences about effects of
causes and about causes of effects, and remarked ``{\ldots}as a
general rule, the effects of causes are far more accessible to our
study than the causes of effects\ldots'' \citep[Book~3, Chapter~10,
\S8]{mill}.  Although a similar distinction has sometimes been
expressed in statistical contexts (see \eg\ \citet{pwh:jasa}), Mill's
associated warning has largely gone unheeded.  We consider that it
deserves more careful attention.  Though evidently related in some
way, problems of CoE are distinct from problems of EoC; indeed, as
Mill understood, they are considerably more subtle and difficult to
handle.

In this article, which builds on and extends \citet{apd:aberdeen}\ and
\citet{fienberg2014}, we attempt to delineate both the differences and
the connexions between these two distinct inferential enterprises.  An
understanding of these issues will clearly be crucial if generic
\adds{retrospective observational} EoC evidence, such as that of the
\citet{frachon:plos1}\ study, is to be brought to bear on an
individual CoE case, such as the toxic tort case of \secref{hypo}. In
particular we shall consider the possibilities of using statistical
evidence to inform CoE inferences.

\subsection{Aspirin trial}
\label{sec:aspirin}

As a simple concrete example, we contrast the following two questions:

\begin{description}
\item[Effects of Causes (EoC)] Ann has a headache.  She is wondering
  whether to take aspirin.  Would that cause her headache to disappear
  (within, say, 30 minutes)?
\item[Causes of Effects (CoE)] Ann had a headache and took aspirin.
  Her headache went away after 30 minutes.  Was that caused by the
  aspirin?
\end{description}
Note that---in a departure from previous related treatments---in both
questions we have separated out the r\^oles of the subject (``Ann''),
on whom we have some information, and the questioner or analyst
(henceforth ``I''), who wants to interpret that information: these
could be the same individual, but need not be.  Any uncertainty about
the answers to the above queries is my personal uncertainty, and is
most properly regarded as a subjective probability, though informed by
relevant data.  This is somewhat analogous to the situation in court,
where we distinguish between a witness, who supplies evidence (\eg, on
epidemiology), and the trier of fact, be it a judge or a jury, who has
to assess the uncertainty to associate with the question of ultimate
legal interest: the cause of the effect.

What might be relevant data in the present instance? We suppose that a
well-conducted (large, prospective, randomised, double-blind,\ldots)
comparative clinical trial has indicated the following recovery rates:
\begin{eqnarray}
  \label{eq:asp1}
  \Pr(R = 1\cd E = 1) & = & 30\%\\
  \label{eq:asp0}
  \Pr(R = 1\cd E = 0) & = & 12\%
\end{eqnarray}
where $E=1$ [resp., $0$] denotes ``exposure to'' ( = treatment with)
aspirin [resp., no aspirin], and $R=1$ [resp., $0$] denotes that the
headache does [resp., does not] disappear (within 30 minutes).  Here
and throughout, we use $\Pr(\cdot)$ to denote probabilities
(henceforth termed {\em chances\/}) underlying a population
data-generating process.

\section{Statistical Evidence for EoC}
\label{sec:stateoc}

\begin{quote}
  % \item[Effects of Causes (EoC)]
  \bf Ann has a headache.  She is wondering whether to take aspirin.
  Would that cause her headache to disappear (within, say, 30
  minutes)?
\end{quote}
Most of classical statistical experimental design and inference is
geared to elucidating the effects of causes, and much careful
attention over many years has gone into clarifying and improving
methods for doing this, for example by the use of randomised
comparative experiments \citep{fisher:doe,hill:clintrial} to control
for potential confounding factors.  Even when emphasis is specifically
targeted on statistical causality
\takep{\citep{dbr:jep,pearl:book}}\addp{\citep{dbr:jep,HollandRubin1988,pearl:book}}
this still mostly addresses EoC problems, albeit in observational
rather than experimental settings.

In order to highlight the major issue, we confine attention here to
data from a study, such as the aspirin trial of \secref{aspirin}, that
can be regarded as supporting genuinely causal
inferences.\footnote{Some considerations relevant to the possibilities
  for causal inference in various data-collection settings can be
  found in \citet{apd:aberdeen}.}  In particular, for the aspirin
trial this would mean that---so long as I can regard Ann as being
comparable with the patients in the trial---if she takes aspirin I can
expect her headache to disappear within 30 minutes with probability
$30\%$, but with probability only $12\%$ if she does not.  If I myself
am Ann, then (other things being equal) taking the aspirin is my
preferred option.

In this case, the EoC causal inference is based on a simple contrast
between the two ``prospective'' conditional probabilities, $\Pr(R =
1\cd E = 1)$ and $\Pr(R = 1\cd E = 0)$. In particular, the information
needed for making EoC causal inferences---and so for guiding future
decisions---is subsumed in the \takem{joint} \addm{conditional}
probability distribution of the \takep{exposure $E$ and the }response
$R$ given \addp{exposure $E$}.  In more complex situations we may have
to make various modifications, \eg\ adjustment for covariates, but the
essential point remains that purely probabilistic knowledge, properly
conditioned on known facts, is sufficient to address EoC-type
questions.

% \takem{ \section{How to Understand ``Causes of Effects''?}}
% \label{sec:coe}

\takem{Addressing a CoE-type question is much more
  problematic---indeed, even to formulate the question clearly is a
  nontrivial enterprise.  We can no longer base our approach purely on
  the probability distribution of $E$ and $R$ conditioned on known
  facts, since we know the values of both variables ($E=1$, $R=1$),
  and after conditioning on that knowledge there is no probabilistic
  uncertainty left to work with.

  One possible approach, popular in statistical circles, is based on
  the concept of the ``counterfactual contrast'', which in turn rests
  on the introduction of ``potential responses'' \citep{dbr:jep}.  We
  proceed by splitting the response variable $R$ into two variables,
  $R_0$ and $R_1$, where we conceive of $R_1$ [resp., $R_0$] as a
  potential value of $R$, that will eventuate if in fact $E=1$ [resp.,
  $0$].  Both these potential responses are regarded as existing prior
  to the determination of $E$. We thus now need to model the three
  variables $(E, R_0, R_1)$ together, rather than (as previously) just
  the two variables $(E, R)$.\footnote{The observed response $R$ is
    determined by these three variables as $R = R_E$.}  We might now
  cast the CoE question as enquiring about the relationship between
  $R_0$ and $R_1$.  Thus ``$R_1 = 1, R_0 = 0$'' describes the
  situation where Ann's headache disappears if she takes the aspirin,
  but does not if she does not---a state of affairs that might
  reasonably be described as the disappearance of Ann's headache being
  {\em caused\/} by taking the aspirin.  In particular, if Ann has
  taken the aspirin and her headache disappeared (thus $R_1 = 1$),
  these two events can be regarded as causally connected just in the
  case that $R_0 = 0$.}

\takem{\subsection{Science and Policy}
  \label{sec:policy}

  Although we shall follow through with the above formulation in the
  remainder of this article, we here turn aside to consider an
  objection to it: it simply might not be appropriate to regard, as
  the ``counterfactual foil'' to the factual response ($R_1$), what
  would have happened ($R_0$) if the exposure had not occurred ($E=0$)
  but all other prior circumstances were the same.  For example, there
  has been a series of legal cases in which various administrations
  %%%% in both the UK and the US?
  have sued tobacco companies on the basis that they had not properly
  informed the public of the dangers of smoking when they first had
  that evidence, and should therefore be liable for the increased
  costs that fell on health services due to that act of omission.  But
  it could be argued that, since smokers tend to die earlier than
  non-smokers, encouraging (or at least not discouraging) smoking
  would in fact reduce the total burden on the health services.  Such
  an attempted defence has, however, usually been ruled inadmissible.
  Instead, as a matter of policy, the relevant counterfactual
  comparator is taken to be a hypothetical universe in which every one
  lives just as long as they do in fact, but they are healthier
  because they smoke less.  Here we see Science and Policy as
  inextricably intertwined in formulating the appropriate CoE
  question.  And the conceptual and implementational difficulties that
  we discuss below, that beset even the simplest case of inference
  about causes of effects, will be hugely magnified when we wish to
  take additional account of such policy considerations.  }

\section{Statistical Evidence for CoE}
\label{sec:coe}
\begin{quote}
  % \item[Causes of Effects (CoE)]
  \bf Ann had a headache and took aspirin.  Her headache went away
  after 30 minutes.  Was that caused by the aspirin?
\end{quote}
\addm{\subsection{How to Understand ``Causes of Effects''?}

  Addressing a CoE-type question is much more problematic---indeed,
  even to formulate the question clearly is a nontrivial enterprise.
  We can no longer base our approach purely on the probability
  distribution of $E$ and $R$ conditioned on known facts, since we
  know the values of both variables ($E=1$, $R=1$), and after
  conditioning on that knowledge there is no probabilistic uncertainty
  left to work with.

  One possible approach, popular in statistical circles, is based on
  the concept of the ``counterfactual contrast'', which in turn rests
  on the introduction of ``potential responses'' \citep{dbr:jep}.  We
  proceed by splitting the response variable $R$ into two variables,
  $R_0$ and $R_1$, where we conceive of $R_1$ [resp., $R_0$] as a
  potential value of $R$, that will eventuate if in fact $E=1$ [resp.,
  $0$].  Both these potential responses are regarded as existing prior
  to the determination of $E$. We thus now need to model the three
  variables $(E, R_0, R_1)$ together, rather than (as previously) just
  the two variables $(E, R)$.\footnote{The observed response $R$ is
    determined by these three variables as $R = R_E$.}

  We might now cast the CoE question as enquiring about the
  relationship between $R_0$ and $R_1$.  Thus ``$R_1 = 1, R_0 = 0$''
  describes the situation where Ann's headache disappears if she takes
  the aspirin, but does not if she does not---a state of affairs that
  might reasonably be described as the disappearance of Ann's headache
  being {\em caused\/} by taking the aspirin.  In particular, if Ann
  has taken the aspirin and her headache disappeared (thus $R_1 = 1$),
  these two events can be regarded as causally connected just in the
  case that $R_0 = 0$.

  \subsection{Science and Policy}
  \label{sec:policy}

  Although we shall follow through with the above formulation in the
  remainder of this article, we here turn aside to consider an
  objection to it: it simply might not be appropriate to regard, as
  the ``counterfactual foil'' to the factual response ($R_1$), what
  would have happened ($R_0$) if the exposure had not occurred ($E=0$)
  but all other prior circumstances were the same.  For example, there
  has been a series of legal cases in which various administrations
  %%%% in both the UK and the US?
  have sued tobacco companies on the basis that they had not properly
  informed the public of the dangers of smoking when they first had
  that evidence, and should therefore be liable for the increased
  costs that fell on health services due to that act of omission.  But
  it could be argued that, since smokers tend to die earlier than
  non-smokers, encouraging (or at least not discouraging) smoking
  would in fact reduce the total burden on the health services.  Such
  an attempted defence has, however, usually been ruled inadmissible.
  Instead, as a matter of policy, the relevant counterfactual
  comparator is taken to be a hypothetical universe in which every one
  lives just as long as they do in fact, but they are healthier
  because they smoke less.  Here we see Science and Policy as
  inextricably intertwined in formulating the appropriate CoE
  question.  And the conceptual and implementational difficulties that
  we discuss below, that beset even the simplest case of inference
  about causes of effects, will be hugely magnified when we wish to
  take additional account of such policy considerations.
  
  \subsection{Statistical Evidence} }
\label{sec:statcoe}

After the above detour, we return to our formulation of the CoE
question, in terms of a contrast between $R_1$, the actually observed
response (in this case, $R_1=1$) to the treatment actually taken
($E=1$), and $R_0$, the (necessarily unknown) counterfactual response,
that would have been observed had Ann in fact not taken the aspirin.
If ``in counterfact'' $R_0 = 1$, then Ann's headache would have
disappeared even if she had not taken the aspirin, so I must conclude
that it was not the aspirin that cured her.  Conversely, if $R_0=0$
then I can indeed attribute her cure to having taken the aspirin.  In
this way, we formulate the CoE causal question in terms of the
contrast between the factual outcome $R_1$ and the counterfactual
outcome $R_0$.

To address the CoE question I thus need to query $R_0$.  Since $R_0$
has not been observed, it retains a degree of uncertainty, which I
could try to express probabilistically.  However, not only have I not
observed $R_0$, there is, now, no way I could ever observe it, since
once I have observed $R_1$, $R_0$ has becomes a counterfactual
quantity, predicated on a condition ($E=0$) that is counter to known
facts ($E=1$).  This logical difficulty leads to a degree of
unavoidable ambiguity affecting our ability to address the CoE
question.

In evaluating my probabilistic uncertainty, I should condition on all
I know.  My full knowledge about Ann can be expressed as $(E=1, R_1
=1, H)$, where $H$ denotes all the background knowledge I have about
Ann, and the other variables are likewise individualised to her.  With
this understanding, we formally define my PROBABILITY OF
CAUSATION\takep{\footnote{This is similar to what
    \citet[Chapter~9]{pearl:book}\ terms the ``Probability of
    Necessity'' (PN)---which however does not account for the
    conditioning on $H$.  See also \citet{tianpearl2000}.}} as the
{\em conditional probability\/}:
\begin{equation}
  \label{eq:pc00}
  \pc_A = \p_A(R_0 = 0 \cd H, E=1, R_1 = 1)
\end{equation}
where $\p_A$ denotes my probability distribution over attributes of
Ann.

But how can I go about evaluating $\pc_A$, and what other evidence
could be used, and how, to inform this evaluation?  In particular,
how---if at all---could I make use of EoC probabilities such as
\eqref{asp1} and \eqref{asp0} to assist my evaluation of the CoE
probability \eqref{pc00}?

\subsection{Bounding the probability of causation}
\label{sec:ass}
We note that \eqref{pc00} involves a joint distribution of $(R_0,
R_1)$.  Since, as a matter of definition, it is never possible to
observe both $R_0$ and $R_1$ on the same individual, it is problematic
to estimate such a joint distribution.  We might however have a hope
of assessing separate marginal probabilities for $R_0$ and $R_1$; and
this information can be used to set bounds on \pc.  Indeed it is
straightforward to show (\cf \citet{apd:aberdeen}): \takep{
  \begin{equation}
    \label{eq:rrA}
    \min\left\{1,  \frac{1}{\p_A(R_1 = 1 \cd H, E=1)} - \frac{1}{\rr_A}\right\}
    \geq \pc_A \geq \max\left\{0,1 - \frac{1}{\rr_A}\right\},
  \end{equation}
} \addp{
  \begin{equation}
    \label{eq:rrA}
    \min\left\{1,  \frac{\p_A(R_0 = 0 \cd H, E=1)}{\p_A(R_1 = 1 \cd H, E=1)}\right\}
    \geq \pc_A \geq \max\left\{0,1 - \frac{1}{\rr_A}\right\},
  \end{equation}
} where
\begin{equation}
  \label{eq:exptrre}
  \rr_A := \frac{\p_A(R_1 = 1 \cd  H, E=1)}{\p_A(R_0 = 1 \cd H, E=1)}.
\end{equation}
% is the {\em causal risk ratio\/}.
Readers will recognize (\ref{eq:rrA}) as a version of the
Bonferroni-Fr\'echet-Hoeffding bounds
\citep{bonferroni1936,frechet1940,hoeffding1940}\ that play important
r\^oles in other areas of statistics \adds{such as in the study of
  copulas}.

The inequality (\ref{eq:rrA}) will yield a non-trivial lower bound so
long as $\rr_A > 1$, which we can interpret as saying that there is a
positive causal effect of exposure on outcome: \cf the related
argument in \citet{robinsgreenland1989}.  Whenever $\rr_A$ exceeds 2,
we can deduce from \eqref{rr}, without making any further assumptions,
that $\pc_A$ must exceed $50\%$.  In a civil legal case such as that
of \secref{hypo}, causality might then be concluded ``on the balance
of probabilities''.  It is however important to note \addm{(see
  \citet{robins1989})} that, when $\rr_A < 2$, it would not be correct
to conclude from this that $\pc_A < 50\%$ (which would lead to the
case failing); rather, we can only say that we can not be sure that
the probability of causation exceeds $50\%$.

The upper bound in \eqref{rrA} is more subtle.  It is less than $1$
when $\p_A(R_0 = 1 \cd H, E=1) + \p_A(R_1 = 1 \cd H, E=1) > 1$.  This
happens in general only when both Ann's potential outcomes \takep{are
  ``highly likely.''}\addp{have a substantial probability of taking
  value $1$.}  If $\p_A(R_1 = 1 \cd H, E=1) = \p_A(R = 1 \cd H, E=1)$
is only modest in size, \eg, less than 1/2 and $\rr_A > 1$, then the
upper bound is 1.  If $\rr_A$ is large, \eg, $\rr_A>10$, the upper
bound will again be 1 unless $\p_A(R = 1 \cd H, E=1)$ is close to 1.
For the remainder of the paper, for simplicity we proceed using an
upper bound of 1.  Thus we work with the bounds
\begin{equation}
  \label{eq:bounds}
  1\geq \pc_A \geq \max\left\{0,1 - \frac{1}{\rr_A}\right\}
\end{equation}
with $\rr_A$ given by \eqref{exptrre}.

\subsection{The risk ratio}
\label{sec:estbounds}

\addp{Expression \eqref{pc00} and the}\takep{The} denominator of
\eqref{exptrre} involve\takep{s} a counterfactual consideration: of
$R_0$, Ann's potential response were she not to have taken the
aspirin, in the situation that she is known to have taken aspirin
($E=1$).  So it would seem problematic to \takep{estimate
  it}\addp{attempt to identify these quantities} from data.  However,
if my background knowledge $H$ of Ann (on which my distribution $\p_A$
is being conditioned) is sufficiently detailed, then, at the point
before Ann has decided whether or not to take the aspirin, it might
seem appropriate to consider that my uncertainty, conditional on $H$,
about the way her treatment decision $E$ will be made would not
further depend on the (so far entirely unobserved) potential responses
$(R_0,R_1)$.  That is, in this case we might assume
\begin{equation}
  \label{eq:ci}
  (R_0,R_1) \cip_A E \cd H
\end{equation}
where $\cip_A$ denotes conditional independence \citep{apd:CIST}\ in
my distribution $\p_A$ for Ann's characteristics.  When \eqref{ci}
holds we will term the background information $H$ {\em sufficient\/}.
\addp{Then \eqref{pc00} becomes
  \begin{equation}
    \label{eq:pc01}
    \pc_A = \p_A(R_0 = 0 \cd H, R_1 = 1)
  \end{equation}
  and in the lower bound in \eqref{bounds} we can replace
  \eqref{exptrre} by
  \begin{equation}
    \label{eq:rr}
    \rr_A = \frac{\p_A(R_1 = 1 \cd H)}{\p_A(R_0 = 1 \cd H)},
  \end{equation}
  my {\em causal risk ratio\/} for Ann.}\footnote{\addp{We can derive
    \eqref{rr}--though not \eqref{pc01}--from the weaker condition
    that replaces the joint property \eqref{ci} by the two marginal
    properties $R_j \cip_A E \cd H$, $j=0,1$.  Since we are only
    concerned with bounds in this paper, that weaker condition would
    be adequate for our purposes.  However, we find it hard to imagine
    circumstances where we would be willing to accept the weaker but
    not the stronger condition, so will continue to use conditions
    like \eqref{ci}.}}

Sufficiency is a kind of ``no confounding'' requirement on my
distribution $\p_A$ for Ann \addm{(see for instance
  \citet{Dawid2014})}.  It would fail if, for example, I thought that
Ann might take the treatment if she felt really poorly, and not
otherwise; but I did not initially have information as to how she
felt.  Then observing that she took the treatment ($E=1$) would inform
me that she was feeling poorly, so decreasing the probability of a
good response (whether actual, $R_1$, or counterfactual, $R_0$).  Now
if I myself am Ann, my $H$ will already include my own knowledge of my
perceived state of health, so this argument does not apply, and
sufficiency is an acceptable condition.  If I am an external observer,
however, the sufficiency condition is much more problematic, since I
must be able to satisfy myself that my knowledge $H$ of Ann is
complete enough to avoid the above possibility of confounding.  If I
can not assume sufficiency, I can not replace the counterfactual
denominator of \eqref{exptrre} by anything even potentially estimable
from data.

Note that the ``no confounding'' property of sufficiency relates
solely to Ann and my knowledge of her.  It should not be confused with
the superficially similar no confounding property of {\em
  exogeneity\/} described in \secref{estrr} below, which refers, not
to Ann, but to the process whereby possibly relevant data on other
individuals have been gathered.

\subsection{Estimating the risk ratio}
\label{sec:estrr}

Henceforth we assume sufficiency, which at least gets us started, and
\adds{we} aim to see what further progress can be made, and under what
conditions, to get a handle on the bounds on $\pc_A$ supplied by
$\rr_A$.  It is important to be explicit about the assumptions
required, which can be very strong and not easy to justify!

It would be valuable if the probabilities featuring in \eqref{exptrre}
could be related in some way to chances such as \eqref{asp1} and
\eqref{asp0} that are estimable from data.  Consider first the
numerator, the Ann-specific probability $\p_A(R_1 = 1 \cd H, E=1) =
\p_A(R = 1\cd H, E = 1)$.  It is tempting to replace this by the
analogous chance, $\Pr(R = 1\cd H, E = 1)$, which could be estimated
from data as for \eqref{asp1}, based on the subset of treated trial
subjects sharing the same $H$-value as Ann.\footnote{Alternatively,
  the estimate might be constructed from a model for the dependence of
  the response $R$ on $H$ and $E=1$, fitted to all the data, and
  applied with Ann's value of $H$.  We might also be able to reduce to
  a smaller information set $H$, if that is all that is relevant for
  prediction of the responses.}\takep{.}  This would be justified if
we could make the following {\em bold assumption\/} (where Bayesians
can replace the intuitive term ``comparable'' with the more precise
term ``exchangeable''):
\begin{cond}
  \label{cond:treated} Conditional on my knowledge of the
  pre-treatment characteristics of Ann and the trial subjects, I
  regard Ann's potential responses as comparable with those of the
  treated subjects having characteristic $H$.
\end{cond}
Up to this point we have not needed the assumption that $H$ is
sufficient.  But consider now the denominator of \eqref{rrA}.  Because
of its counterfactual nature, we can not argue directly as above.
However, with sufficiency of $H$ we have $\p_A(R_0 = 1 \cd H, E=1) =
\p_A(R_0 = 1\cd H, E=0) = \p_A(R = 1\cd H, E=0)$; and we can estimate
this from the \adds{clinical} trial data, \eg as the estimated chance
${\Pr}(R = 1\cd H, E=0)$, if we can assume:
\begin{cond}
  \label{cond:untreated} Conditional on my knowledge of the
  pre-treatment characteristics of Ann and the trial subjects, I
  regard Ann's potential responses as comparable with those ot the
  untreated subjects having characteristic $H$.
\end{cond}

Now if both \condref{treated} and \condref{untreated} are to hold,
then (by Euclid's first axiom, ``Two things that are equal to the same
thing are also equal to each other''), the groups of trial subjects
with Ann's characteristics $H$ in both arms must be comparable with
each other.  This requires that $H$ be {\em exogenous\/}, in the sense
that, conditional on $H$, the potential outcomes $(R_0,R_1)$ have the
same distribution among treated and untreated study subjects.  This
will hold for a suitably randomised study, and also in certain
observational studies where the possibility of further confounding
factors can be discounted.

Note however that we can not take, as $H$, just {\em any\/} exogenous
set of variables.  The full set of required conditions is:
\begin{enumerate}
\item \label{it:c1} $H$ is exogenous.
\item \label{it:c2} $H$ is sufficient for Ann's response.
\item \label{it:c3} Conditional on $H$, Ann's potential responses are
  comparable with those of the trial subjects.
\end{enumerate}
\addp{We will refer to this set of conditions as the {\em fundamental
    conditions\/}.}

\takep{Only when}\addp{When} we can make good arguments for the
acceptability of \takep{all these strong conditions can we justify
  estimating $\rr_A$ by the population counterpart of \eqref{rr}, the
  {\em observational risk ratio\/}:
  \begin{equation}
    \label{eq:obsrr}
    \orr = \frac{\Pr(R = 1 \cd H, E=1)}{\Pr(R = 1 \cd H, E=0)}.
  \end{equation}}%
\addp{these fundamental conditions, equation~\eqref{pc00} becomes
  \begin{equation}
    \label{eq:pc000}
    \pc_A = \Pr(R_0 = 0 \cd H, R_1 = 1),
  \end{equation}
  and, in the lower bound in \eqref{rrA}, we can identify $\rr_A$ with
  the population counterpart of \eqref{rr}, the {\em observational
    risk ratio\/}:
  \begin{equation}
    \label{eq:obsrr}
    \orr: = \frac{\Pr(R = 1 \cd H, E=1)}{\Pr(R = 1 \cd H, E=0)}.
  \end{equation}

  Now that we have made clear that the fundamental conditions can be
  expected to hold only in special circumstances, when they will
  require detailed justification, we shall henceforth confine
  ourselves to futher consideration of just these special cases.  In
  particular we shall accept \eqref{pc000}, and $\rr_A = \orr$ as in
  \eqref{obsrr}.  So we will use the bounds
  \begin{equation}
    \label{eq:ineqpc}
    1\geq \pc \geq \max\left\{0,1 - \frac{1}{\orr}\right\}
  \end{equation}
  with $\orr$ given by \eqref{obsrr}.  \addp{(Here and henceforth,
    unless the context requires otherwise we drop the identifier $A$
    on $\pc$: these bounds will apply to any individual for whom the
    fundamental conditions hold.)}  }

\addp{

  \subsection{An alternative approach}
  \label{sec:pearl}
  Our probability of causation, $\pc_A$ given by \eqref{pc00}, is
  essentially the same as what \citet[Chapter~9]{pearl:book}\ terms
  the ``Probability of Necessity'' (PN).  \citet{tianpearl2000} take
  an alternative approach to supplying bounds for PN, based on data
  and assumptions different from ours.  In particular, they drop our
  requirement that $H$ be sufficient for Ann's response, requiring
  instead the availability of two sets of data on individuals
  comparable to Ann: one set in which treatment was (or can be
  regarded as) randomized, and another in which it arose ``naturally''
  in the same way as for Ann.  Because of these differences it is not
  in general possible to compare their bounds and ours.  See
  \citet{pearlcomm:smr,dff:pearlresp} for further discussion of these
  issues.

}

\subsection{Uncertain exposure}
\label{sec:uncexp}

\takep{Above}\addp{So far} we have supposed we know both the fact of
exposure ($E=1$) and the fact of response ($R=1$), the only
uncertainty being about whether there was a causal link between these
two facts.  There are other situations where we might observe the
response, and wonder whether it was caused by exposure, without
knowing with certainty whether or not that exposure had in fact taken
place.  In such cases we have to multiply the probability of causation
$\pc_A$ by the probability of exposure, conditional on the known fact
of a positive response, yielding a modified probability of causation:
\begin{equation}
  \label{eq:pc*}
  \pc_A^* = \pc_A \times \p_A(E=1 \cd H, R=1).
\end{equation}

In particular, \takep{when the strong conditions of \secref{estrr},
  justifying the use of population chances $\Pr(\cdot)$ in place of
  Ann-specific probabilities $\p_A(\cdot)$, can be assumed}\addp{under
  the fundamental conditions}, combining this with \takep{\eqref{pc*}
  (and using the upper bound $1$ of
  \eqref{bounds})}\addp{\eqref{ineqpc}} delivers the inequalities
\begin{equation}
  \label{eq:ineqpc*}
  \Pr(E = 1 \cd H, R=1) \geq \pc^* \geq \max\left\{0,1 - \frac{\Pr(E = 0 \cd H, R = 1)}{\Pr( E=0\cd H)}\right\}.
\end{equation}
\takep{(Here and henceforth we drop the identifier $A$ on $\pc^*$:
  these bounds will apply to any individual for which the required
  conditions hold.)}

\section{Statistical Uncertainty}
\label{sec:statunc} Our discussion so far has treated estimates, such
as those in \eqref{asp1} and \eqref{asp0}, as if they were the true
values of the chances.  Even so, we found that we obtain, at best,
only partial CoE information, which confines $\pc$ or $\pc^*$ to an
interval but does not yield a point value.  In real applications our
data will not be extensive enough to give us pinpoint estimates of
even the bounds featuring in the\takep{se} inequality formulae\addp{
  \eqref{ineqpc} or \eqref{ineqpc*}}, and so we have to take
additional account of the resulting statistical uncertainty.  The
result of our inference is thus an {\em uncertain interval\/} within
which a {\em probability\/} ($\pc$ or $\pc^*$) must lie---thus
compounding three different kinds of uncertainty.  \takep{This is a
  novel form of inferential output, and it is far from clear how best
  to express and display it, and what use to make of it.}

\notep{Here --- or somewhere --- I want to add some material relating
  to the alternative strategy of simply replacing probabilities
  featuring in the bounds by their posterior expectation, as discussed
  in my emails of around 20 August}
  
\notes{I have seen people do this and then ignore the impact and so I
  am interested in seeing the details of how far you would take this.}

Statistical uncertainty, at least, is well studied, and can be
expressed and understood in a variety of different ways, as touted and
debated by the various competing schools of statistical inference.
\addp{The generic problem of inference for a quantity (like $\pc$)
  that, being only partly identified by the data, is subject to
  interval bounds, has been treated from both a classical perspective}
\adds{\citep{manski:2003,manski:2007,stijn2006} and a Bayesian
  perspective
  \citep{greenland:2005,greenland2009,gustafson2005,gustafson2009},
  but these approaches usually involve adding assumptions or data via
  model expansion and \takep{thus }are not directly applicable here.}

We \takep{consider it most straightforward here to}\addp{here} take a
Bayesian approach, \takep{which delivers }\addp{to derive} a joint
\addp{posterior} probability distribution (which, following the
helpful terminology of \citet{best}, we henceforth term a {\em
  credence distribution\/}) for \takep{all the}\addp{the estimable}
unknown chances in the problem\takep{, conditional on observed data}.

One possible \addp{Bayesian} tactic would be to \takep{work with a
  joint credence}\addp{assign a prior} distribution to the
\addp{multivariate parameter, $\phi$ say, comprising the} chances
assigned to the four configurations of $(R_0,R_1)$\addp{ conditioned
  on $H$.  Under the fundamental conditions, $\pc_A$ is a function of
  $\phi$, given by \eqref{pc000}, so that a Bayesian analysis, based
  on such a prior, would deliver a fully determined posterior
  distribution for $\pc_A$.  However,} \takep{While this would deliver
  a seemingly comprehensible inference, in the form of a posterior
  credence distribution for $\pc^*$,} this is problematic: because
$R_0$ and $R_1$ are never simultaneously observable, these joint
chances can not be consistently estimated from data, so that this
``inference'' remains highly sensitive to the specific prior
assumptions made, however extensive the data.  \addp{Alternatively
  put, the parameter $\phi$ describing the joint distribution of
  $(R_0,R_1)$ (given $H$) is not fully identifiable from data; at
  best, only the parameter, $\lambda$ say (a non-invertible function
  of $\phi$), determining the associated marginal distributions of
  $R_0$ and $R_1$, is identifiable.  Then $\lambda$ is a {\em
    sufficient parameter\/} \citep{barankin60}.  For extensive data
  (and a non-dogmatic prior), the posterior distribution of $\lambda$
  will converge to a point mass at its true value, but the posterior
  conditional distribution of $\phi$ given $\lambda$ will be exactly
  the same in the posterior as in the prior
  \citep{kadane1974,apd:CIST}.  In particular the marginal posterior
  distribution of any non-identifiable function of $\phi$ will be
  non-degenerate, and highly dependent on the form of the conditional
  prior for $\phi$ given $\lambda$
  \citep{gustafson2005,gustafson2009,gustafson2012}.}

\takep{Instead}\addp{For these reasons}, we prefer to \takep{work
  with}\addp{assign} a joint credence distribution for the (estimable)
marginal chances \addp{alone}: given sufficient data, of sufficiently
good quality, these will be well estimated and insensitive to prior
assumptions.  The price of this increased statistical precision,
however, is logical imprecision, since from \addp{even perfect
  knowledge of} these chances we can at best derive \addp{interval}
inequalities for $\pc$ or $\pc^*$.  Thus our inference has the form of
a {\em random interval\/} asserted to contain $\pc$ or $\pc^*$.

\addp{

  \subsection{Group or individual inference?}
  \label{sec:altview}

  In the above approach, we considered the probabilities featuring in
  the inequalities \eqref{ineqpc} and \eqref{ineqpc*} as ``objective
  chances'', which we might interpret as limiting relative frequencies
  computed in an appropriate groups of exchangeable individuals.  We
  focused on the \takep{non-degenerate} posterior joint credence
  distribution of these chances, given the available statistical data
  $D$---thus giving rise to a random uncertainty interval for the
  probability of causation, itself regarded as an objective chance.
  We refer to this as the ``group-focused'' approach.

  Another approach to using data to inform the inference about $\pc$
  or $\pc^*$ is to regard these concepts, and the probabilities
  featuring in the bounds for them, themselves as credences,
  quantifying numerically the relevant uncertainty about attributes of
  the specific individual, Ann, on whom we are focusing.  This is the
  ``individual-focused'' approach.  For an interchange on these issues
  and ``group-to-individual'' inference, see the discussion by Dawid
  and the authors' rejoinder in \citet{best}.

  In the individual-focused formulation, the term $\Pr(E = 0 \cd H, R
  = 1)$ in \eqref{ineqpc*}, for example, would be replaced by
  $\p_A(E_A = 0 \cd H_A, R_A = 1, D)$, where the suffix $A$ refers to
  attributes of Ann.  We \addp{now} obtain a non-random uncertainty
  interval for $\pc_A$ (or $\pc_A^*$)---but one that is computed in
  the light of the available evidence, and would be likely to change
  were further data to become available.

  To continue with this example, let $\psi$ denote the chance $\Pr(E=0
  \cd H, R=1)$.  If the individuals in $D$ can be regarded as
  exchangeable with Ann, and we interpret $\psi$ as a limiting
  relative frequency in this exchangeable setting, we will have:
  \begin{eqnarray}
    \nonumber
    \p_A(E_A = 0 \cd H_A, R_A = 1, D) &=& \E\left\{\p_A(E_A = 0 \cd H_A, R_A = 1, \psi) \cd 
      H_A, R_A = 1, D\right\}\\
    \label{eq:psi}
    &=&\E\left(\psi \cd H_A, R_A = 1, D\right).
  \end{eqnarray}
  Often, given the data $D$, the further conditioning in \eqref{psi}
  on the Ann-specific information $(H_A, R_A=1)$ will have negligible
  effect---in which case the desired Ann-specific probability
  $\p_A(E_A = 0 \cd H_A, R_A = 1, D)$ can be approximated by the
  posterior expectation (\ie, conditioned on $D$ alone) of the
  conditional chance $\psi = \Pr(E=0 \cd H, R=1)$.  A similar argument
  applies to any other required credences in the problem.

}

\subsection{Additional issues}
\label{sec:issues}

\notep{I think this section needs some substantial rethinking and
  rewriting, in the light of other things we have said}

\takep{The}\addp{All our} above analysis is predicated on the causal
relevance of the epidemiological data, assuming that we can use the
study to obtain a sound estimate of the causal risk ratio \takep{\rr\
}\addp{$\rr_A$} that features in \eqref{exptrre}.  For example, in a
simple fully randomised study we could use \orr, as given by
\eqref{obsrr}, as a proxy for $\rr$.  But such studies are the
exception in epidemiology, so that the issues in real world settings
where interest is focused on the causes of effects are typically much
more complex.  Thus in the benfluorex example of \secref{frachon},
using the frequencies in \tabref{twobytwo} for this purpose, by
plugging them into the formula for \orr\ and interpreting this as \rr,
would be totally misleading, even if we attempted to account for
statistical uncertainty as described above.  Indeed, as \takes{the
  discussion} \adds{we noted} in
Section~\ref{sec:frachon}\takes{noted}, there are additional problems
in this case: because the study of \citet{frachon:plos1} was
retrospective, and the frequencies in \tabref{twobytwo} could not
\takes{even} be used to estimate \orr, even in the absence of
confounding.  And this problem remains when, admitting the likely
existence of confounding, we conduct a more sophisticated
analysis---such as the multiple logistic regression that produced the
adjusted odds ratio---to try and account for it.  Even when this ploy
can be regarded as successful, still the best we can ever do with
retrospective data is to estimate the causal odds ratio---which will
approximate the desired causal risk ratio, as required for setting the
lower bound on $\pc$, only when the outcome is rare.

The judge in the hypothetical case we pose should therefore be doubly
wary of the relevance of the epidemiological evidence when trying to
assess whether the drug caused the plaintiff's heart disease.

There are even more complex situations where the data are \takep{of
  retrospective sort}\addp{retrospective} and where there are multiple
outcomes of interest and multiple time points for their assessment.  A
notable example comes from \takep{a}\addp{the} continuing effort in
the United States to examine the long-term health effects of exposure
to Agent Orange among US Vietnam veterans.  From 1962 to 1971, the US
military sprayed herbicides over Vietnam.  In 1991 the US Congress
passed the Agent Orange Act, requiring a comprehensive evaluation of
scientific and medical information regarding the health effects of
exposure to Agent Orange and other herbicides used in Vietnam:
\citet{national2011Veterans} is the eighth biennial update
implementing this Congressional mandate.  The report examines
epidemiological studies of the health status of veterans considering a
multiplicity of deleterious effects, \eg, different forms of cancer
and early-onset peripheral neuropathy, and with limited information on
exposure, both at the aggregate and individual level.  A standard tool
in the studies incorporated into this regularly-updated assessment is
the use of adjusted odds-ratios from retrospective logistic regression
analyses.  Identification of a substantial \rr triggers compensation
to veterans for health and disability outcomes associated with
putative exposure.

\section{Case Study}
\label{sec:best} We illustrate our analysis with an example taken from
\citet{best}. The motivating real life case was the diagnosis of abuse
in an infant child, $c$, presenting with an acute life threatening
event (``\alte'')\takep{ and nosebleed (``\bleed'')}.  So now we take
exposure, $E=1$, to denote \abuse, and response, $R=1$ to denote
\takep{the combination of \alte and \bleed}\addp{\alte}.

\subsection{Three tasks}
\label{sec:3tasks}

We can distinguish three tasks that we might wish to address
probabilistically concerning the relationship between exposure and
response in this individual case; these are quite distinct and should
not be confused---although there are of course relationships (far from
trivial) between them.
\begin{description}
\item[Forecasting] If the child is abused, what is the probability the
  child will suffer \alte \takep{ and nosebleed?}  \takep{ ---
    $\p_c(\alte\ \&\ \bleed \cd \abuse)$}\addp{ --- $\p_c(\alte \cd
    \abuse)$}?
\item[Backcasting] If the child suffers \alte\takep{ and nosebleed},
  what is the probability the child was abused?  \takep{ ---
    $\p_c(\abuse \cd \alte\ \&\ \bleed)$}\addp{ --- $\p_c(\abuse \cd
    \alte)$}?
\item[Attribution] If the child suffers \alte\takep{ and nosebleed},
  what is the probability \takep{these were}\addp{this was} caused by
  abuse?
\end{description}
In the above, we have used $\p_c$ to indicate my probabilities for
this child (implicitly conditioned on the background information $H$ I
have about the child).  Even so, \takep{can take}\addp{we have a
  choice between taking} a group-focused approach, in which $\p_c$ is
interpreted as an uncertain chance, \addp{relevant to a group of
  individuals of which this child is one;} or an individual-focused
approach, with $\p_c$ \takep{indicating a credence}\addp{denoting my
  credence about this specific child}.  We start by taking the
group-focused approach: the individual-focused approach will be
considered in \secref{indiv}.

\subsection{Attribution analysis}
\label{sec:attrib}

\citet{best}\ focused on the backcasting task: of assessing whether or
not abuse has in fact taken place, based on the data on the individual
case and on relevant statistical studies.  Their analysis directly
addressed the main substantive concern, since it was the occurrence of
abuse---whether or not it in fact caused the observed signs---that was
at issue.  They did not need to enquire whether or not the observed
signs were {\em caused\/} by abuse.  That attribution question however
will be our focus here.  We note that, since the very fact of abuse is
itself uncertain, we also need to consider the backcasting issue.
This is done by taking, as the relevant probabilistic target of our
inference, the modified probability of causation $\pc^*$, as given by
\eqref{pc*}.

We have described in \secref{statcoe} the many very strong assumptions
that have to be made in order to justify using data to estimate even
the weak interval bounds of \eqref{ineqpc*} for $\pc^*$.  In the
present example, the data used by \citet{best}\ were gleaned from a
search for relevant published studies.  Those identified were of
varying design and quality, and the data extracted from them can in no
sense be regarded as supporting genuine causal inferences --- indeed,
it is not easy to find real examples where the conditions supporting
causal inference of this type could be regarded as satisfied.
\takep{Nevertheless, purely}\addp{Purely} for illustration we shall
proceed as if they are, so that we can use the inequalities of
\eqref{ineqpc*}.  As a further\takep{ highly
  unrealistic}\addp{---admittedly highly unrealistic---}assumption, we
take the sufficient information $H$ to be trivial.  \addp{All these
  imperfections in the data, and in our understanding of the context,
  mean that our analysis must not be taken as delivering a credible
  conclusion in this particular application; however, we hope that, by
  following it through in detail, we may help to clarify the points to
  which attention should given when analysing any similar problem.}

Using a Gibbs sampler implemented in the
\winbugs\hspace{-.3em}$^\copyright$\ software \citep{bugsbook},
\citet{best}\ find the posterior credence distributions for various
conditional chances, based on the data.\footnote{\citet{best}\ conduct
  several alternative analyses, with some of the less reliable data
  values being either included or excluded.  Our own analysis is based
  on the predictive model and data in the combined \winbugs\ code of
  Appendices~B and D of their paper, as for their own Table~4.  This
  analysis targets a case-specific chance, having greater relevance,
  but also more uncertainty, than the overall population-based
  chance.} In particular, they obtain the posterior credence
distribution for the conditional chance $\Pr(E=1 \cd R=1)$, and thus
for $\Pr(E=0 \cd R=1) = 1 - \Pr(E=1 \cd R=1)$, as needed on both sides
of \eqref{ineqpc*}.

For our purposes, however, we need more: the lower bound for $\pc^*$
in \eqref{ineqpc*} also involves the marginal prior chance $\Pr(E =
0)$---or, essentially equivalently, $\theta := \Pr(E = 1)$, the chance
of abuse having taken place (in this individual case), before the
evidence of \alte\ \takep{and \bleed\ }is taken into consideration.
And there is no available statistical evidence relevant to this
quantity.

\begin{figure}[htbp]
  \centering \hfill \subfigure[Prior~1: $\theta \sim\beta(0.1,0.1)$]
  {\includegraphics[scale=.3]{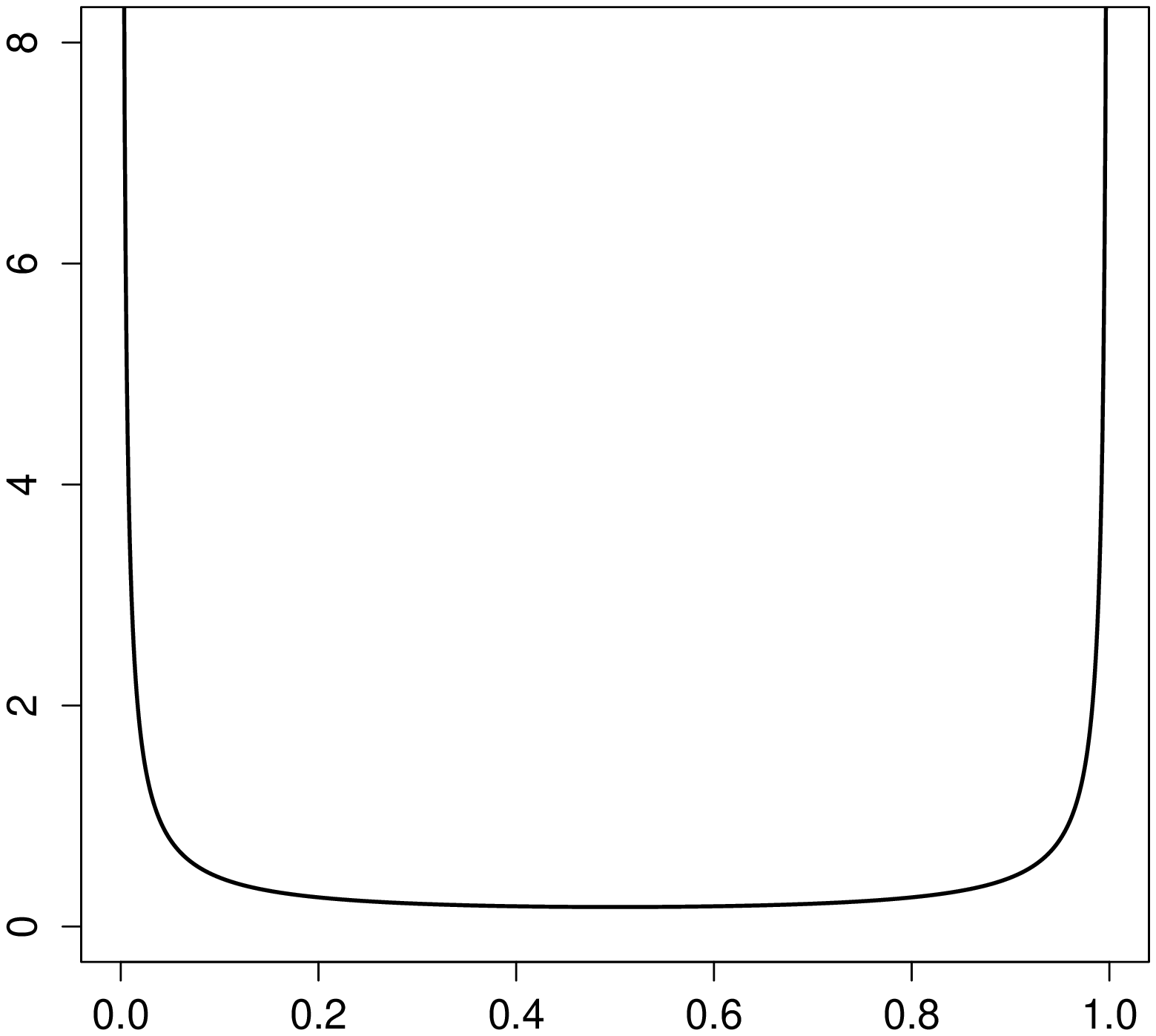}}\hfill \subfigure[Prior~2:
  $\theta
  \sim\beta(1,9)$]{\includegraphics[scale=.3]{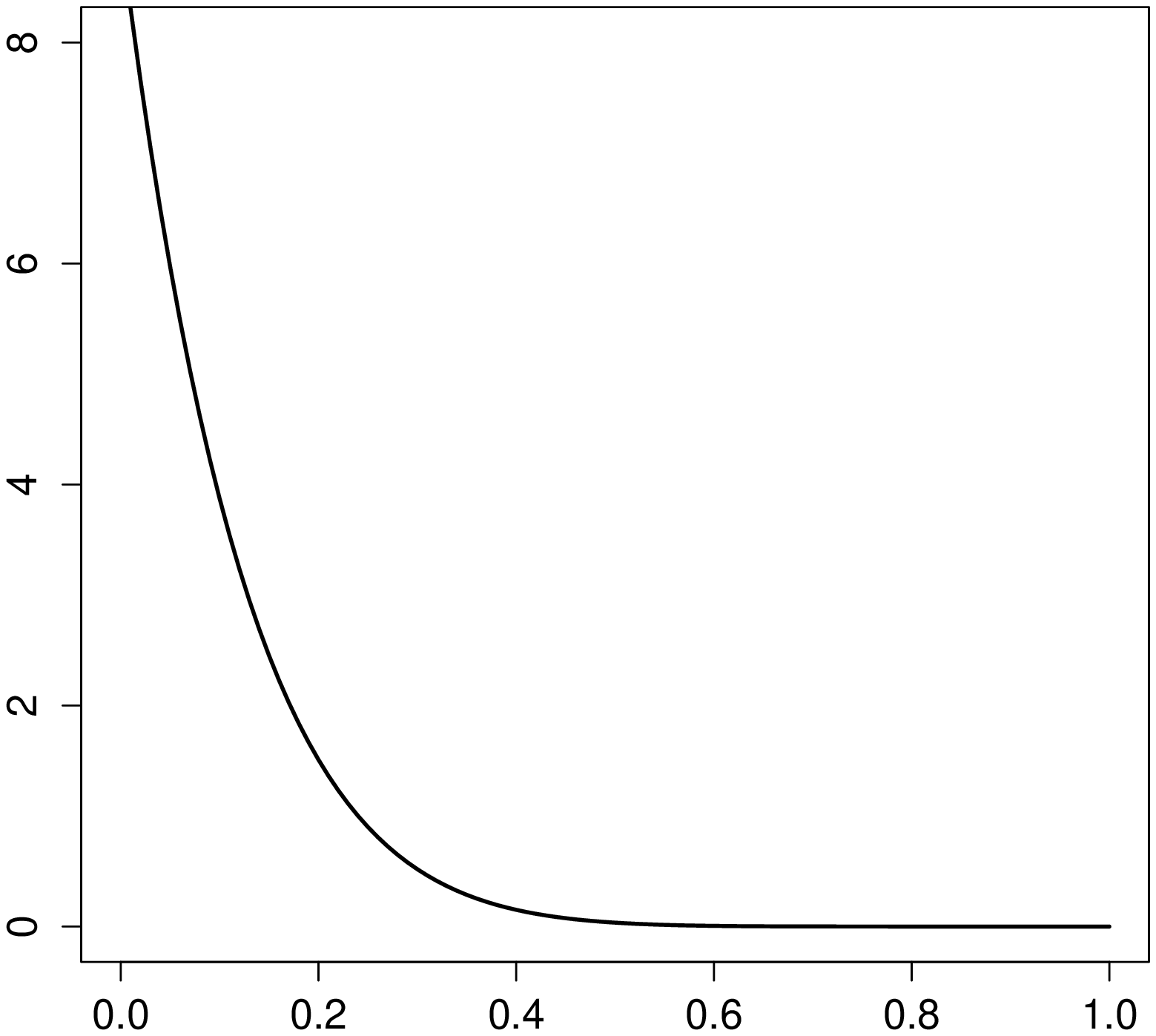}}\hfill\hfill
  \caption{Two prior credence distributions for $\theta$\addp{, the
      prior chance that abuse has taken place.}}
  \label{fig:priors}
\end{figure}

We therefore proceed by introducing our own prior credence
distribution for $\theta$, and treating this chance as independent of
all the others in the problem.  We can expect considerable sensitivity
to the specific choice made.  To begin to explore this, we try two
different prior credence distributions for $\theta$, both beta
distributions for simplicity and tractability:
\begin{description}
\item[Prior~1:] $\theta \sim \beta(0.1, 0.1)$.\\
  This has mean $0.5$ and standard deviation $0.46$.  It can be
  regarded as representing very substantial prior uncertainty about
  $\theta$.
\item[Prior~2:] $\theta\sim \beta(1,9)$.\\
  This has mean $0.1$ and standard deviation $0.09$.  While still
  admitting uncertainty, it attempts to take into account the prior
  unlikelihood of abuse: its mean $0.1$ is the unconditional
  probability assigned to this event.
\end{description}
Density functions of these two prior credence distributions are
displayed in \figref{priors}.

\addp{We note that the lower bound in \eqref{ineqpc*} is $0$ if and
  only if $\rr_A\leq 1$, which event is independent of $\theta$.
  Consequently the posterior probability that the lower bound is $0$
  is unaffected by the assumptions made about the prior distribution
  for $\theta$---as similarly is the conditional posterior
  distribution of the upper bound, given that the lower bound is $0$.}

\section{Data Analysis}
\label{sec:datanal}

We have conducted our own analysis of the data, based on the \winbugs\
code of \citet{best}\ elaborated so as to incorporate $\theta$.
\takem{A chain of length 500000 was generated after a burn-in phase of
  500000 iterations.  Interpreting the generated chain as an
  independent sample of size 500000 would give potentially misleading
  density estimates: for this reason we have based the estimates on
  thinned samples taking every 10th elements of the chain, as
  suggested by an analysis (not reported) of autocorrelation
  estimates.}  \addm{After a burn-in phase of 500000 iterations, to
  get rid of autocorrelation we have based the estimates on thinned
  samples taking every 10th elements of the chain. Then we have
  considered a chain of length 50000.}

\subsection{Bivariate distribution}
\label{sec:biv}

A complete inference would describe the posterior credence
distribution of the \addp{uncertainty} interval \eqref{ineqpc*} for
$\pc^*$, whose end-points are functions of random chances, and hence
themselves have a bivariate distribution.

\begin{figure}[htbp]
  \centering \hfill \subfigure[Prior~1]
  {\includegraphics[scale=.3]{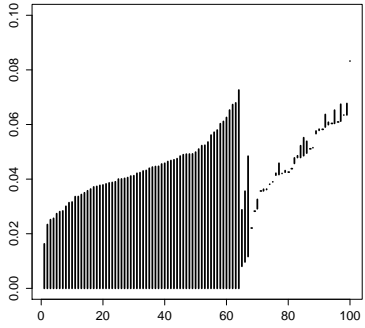}}\hfill \subfigure[Prior~2]
  {\includegraphics[scale=.3]{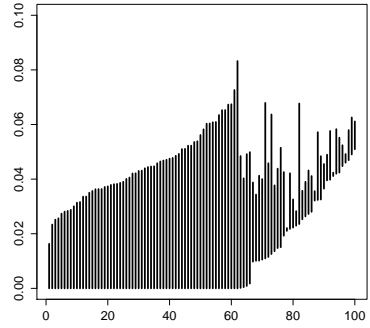}}\hfill\hfill
  \caption{\addp{For each of the priors of \figref{priors}, 100
      uncertainty intervals, randomly sampled from the bivariate
      posterior distribution of the lower and upper bounds, are
      displayed.  The i}\takep{I}ntervals are ordered in increasing
    value of the lower bound.}
  \label{fig:lengths}
\end{figure}

Note that, whenever the inequality $\Pr(E = 1 \cd R = 1) \leq \Pr(
E=1)$ between chances holds, which corresponds to negative association
between exposure and outcome and will happen with positive probability
in the posterior credence distribution, the lower bound of the
\addp{uncertainty} interval is 0 and is thus entirely uninformative.
Thus the posterior credence distribution is a mixture of a continuous
bivariate distribution, and (with positive probability) a distribution
for the upper bound alone.  The probability that the lower bound of
the \addp{uncertainty} interval is $0$ \addp{(which is independent of
  the prior distribution used)} is estimated as \takep{$0.627$ for
  Prior~1, and $0.677$ for Prior~2.}\addp{$0.65$.}  \figref{lengths}
displays, for the two different priors, samples from the bivariate
posterior credence distribution (ordered by lower bound).  In the
plots are reported 100 \addp{uncertainty} intervals obtained by
selecting one iteration of the chain every \takem{5000} \addm{500}.

% A complete inference would be to describe the posterior credence
% distribution of the interval \eqref{ineqpc*} for $\pc^*$, whose
% end-points are functions of random chances, and hence themselves
% have a bivariate posterior credence distribution. It is not easy,
% however, to describe and display this complete inference
% meaningfully. Instead, we describe some summary aspects of it,
% realising that these give only incomplete and inadequate information
% regarding the full inference. {\color{blue} After a burn-in phase of
% 500,000 iterations, we used 500,000 subsequent iterations in our
% computations.  To reduce serial correlation, we selected one
% iteration every 500.}  \todo[inline]{Monica to include various new
% plots here}

\begin{figure}[htbp]
  \centering \hfill \subfigure[Prior~1]
  {\includegraphics[scale=.3]{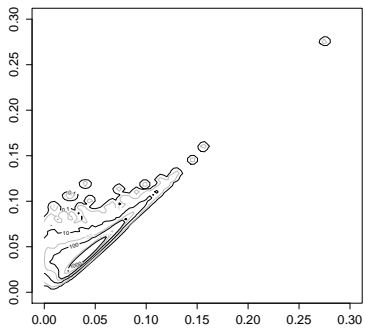}}\hfill \subfigure[Prior~2]
  {\includegraphics[scale=.3]{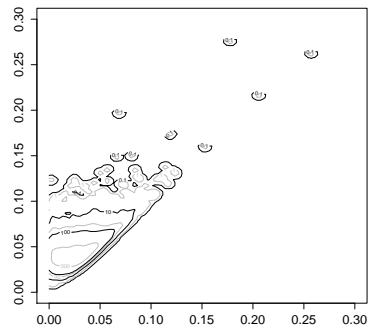}}\hfill\hfill
  \caption{\addp{For each of the priors of \figref{priors}, a }contour
    plot of the joint posterior distribution of \addp{the} lower and
    upper bounds \takep{(for lower bound $>0$)}\addp{of the random
      uncertainty interval, conditional on the lower bound being
      positive.}}
  \label{fig:contour}
\end{figure}

\begin{figure}[htbp]
  \centering %\hfill
  {\includegraphics[scale=.3]{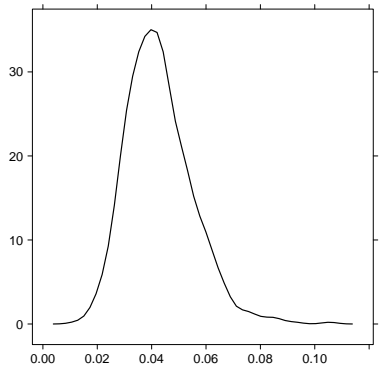}}%\hfill
  \caption{\takep{Posterior density of upper bound (for lower bound
      $=0$}\addp{The posterior density of the upper bound of the
      random uncertainty interval, conditional on the lower bound
      being $0$.  This distribution is independent of the chosen prior
      for the parameter $\theta$.}}
  \label{fig:UB0}
\end{figure}

In \figref{contour} are shown bivariate contour plots, for Priors~1
and 2, of the end-points of the random \addp{uncertainty} interval,
excluding those cases where the lower bound is equal to zero.  The
full joint distribution is completed by specifying the distribution of
the upper bound for these cases\takep{: these are shown in
  \figref{UB0}.}\addp{.  This distribution, which is independent of
  the assumed prior for $\theta$, is shown in \figref{UB0}, which can
  also be interpreted as displaying the conditional distribution of
  the length of the \addp{uncertainty} interval for $\pc^*$, given
  that its lower bound is $0$.}

% As expected prior choice has a big impact on the analysis.

\subsection{Univariate summaries}
\label{sec:univ}

Useful univariate summaries of the overall bivariate inference are the
marginal posterior credence distributions of the upper and lower
bounds, and of the length of the \addp{uncertainty} interval.

\subsubsection{Upper bound}
\label{sec:upper}

\begin{figure}[htbp]
  \centering
  \includegraphics[scale=0.3]{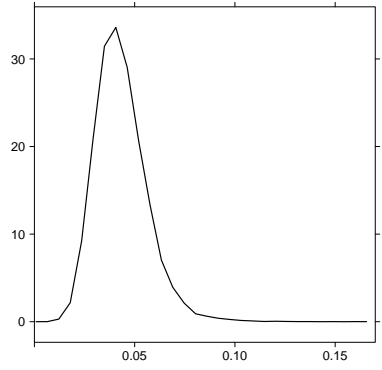}
  \caption{\addp{The p}\takep{P}osterior credence density of the upper
    bound for $\pc^*$.  \addp{This distribution is independent of the
      chosen prior for the parameter $\theta$.}}
  \label{fig:upperdens}
\end{figure}

The upper bound $\Pr(E = 1 \cd R=1)$ in \eqref{ineqpc*} is the chance
of abuse given the case evidence \addp{of \alte}, as already
considered by \citet{best}.  Its posterior credence distribution
(which is unaffected by the choice of prior for $\theta$) is
summarised in the \takep{first}\addp{second} row of Table~4 of
\citet{best}.  We compute the posterior mean and standard deviation
for this upper bound to be $0.043$ and $0.013$, respectively.  Its
posterior density is shown in \figref{upperdens}.

\subsubsection{Lower bound}
\label{sec:lower}

\begin{figure}[htbp]
  \centering \hfill \subfigure[Prior~1]
  {\includegraphics[scale=.3]{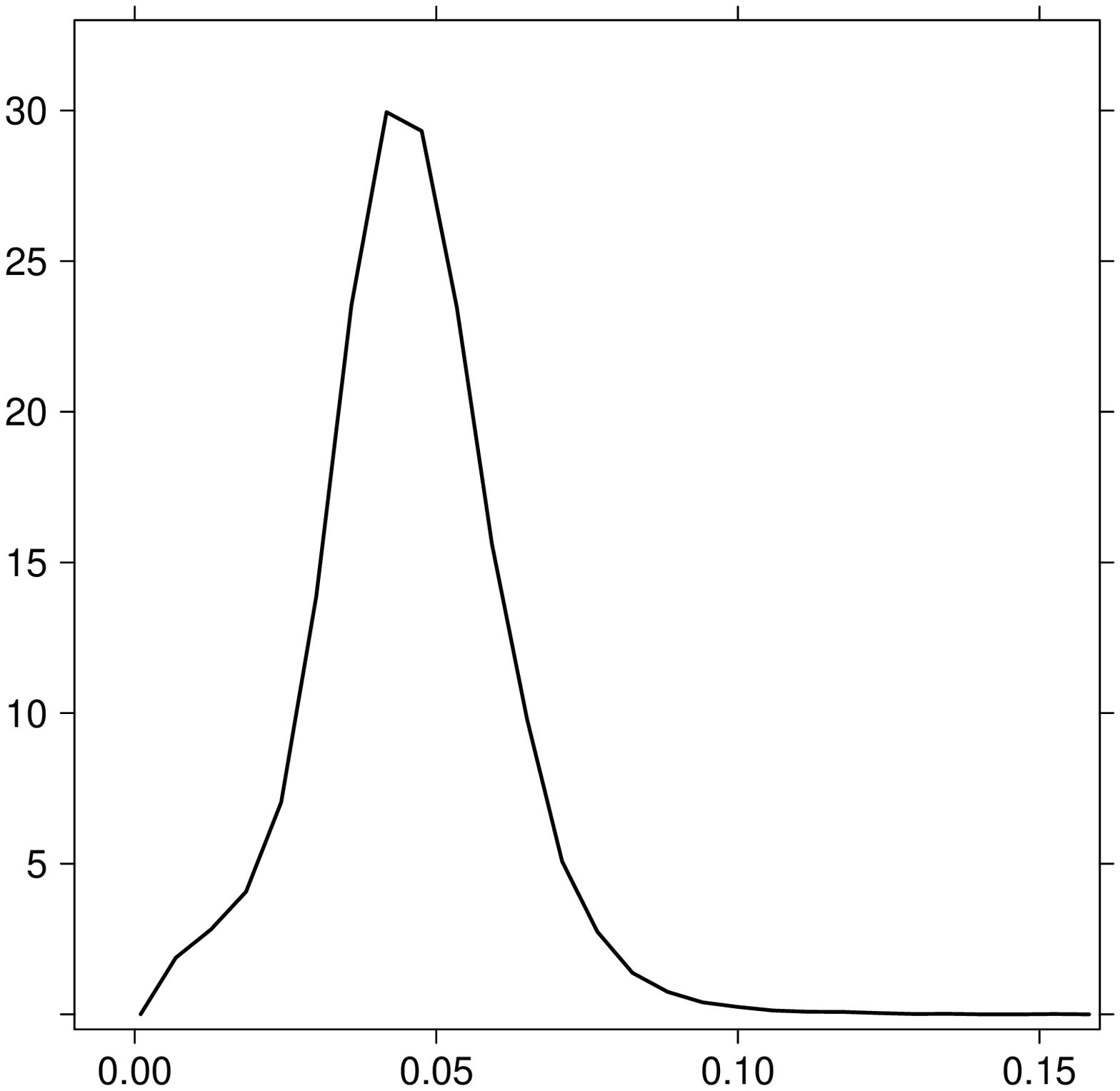}}\hfill \subfigure[Prior~2]
  {\includegraphics[scale=.3]{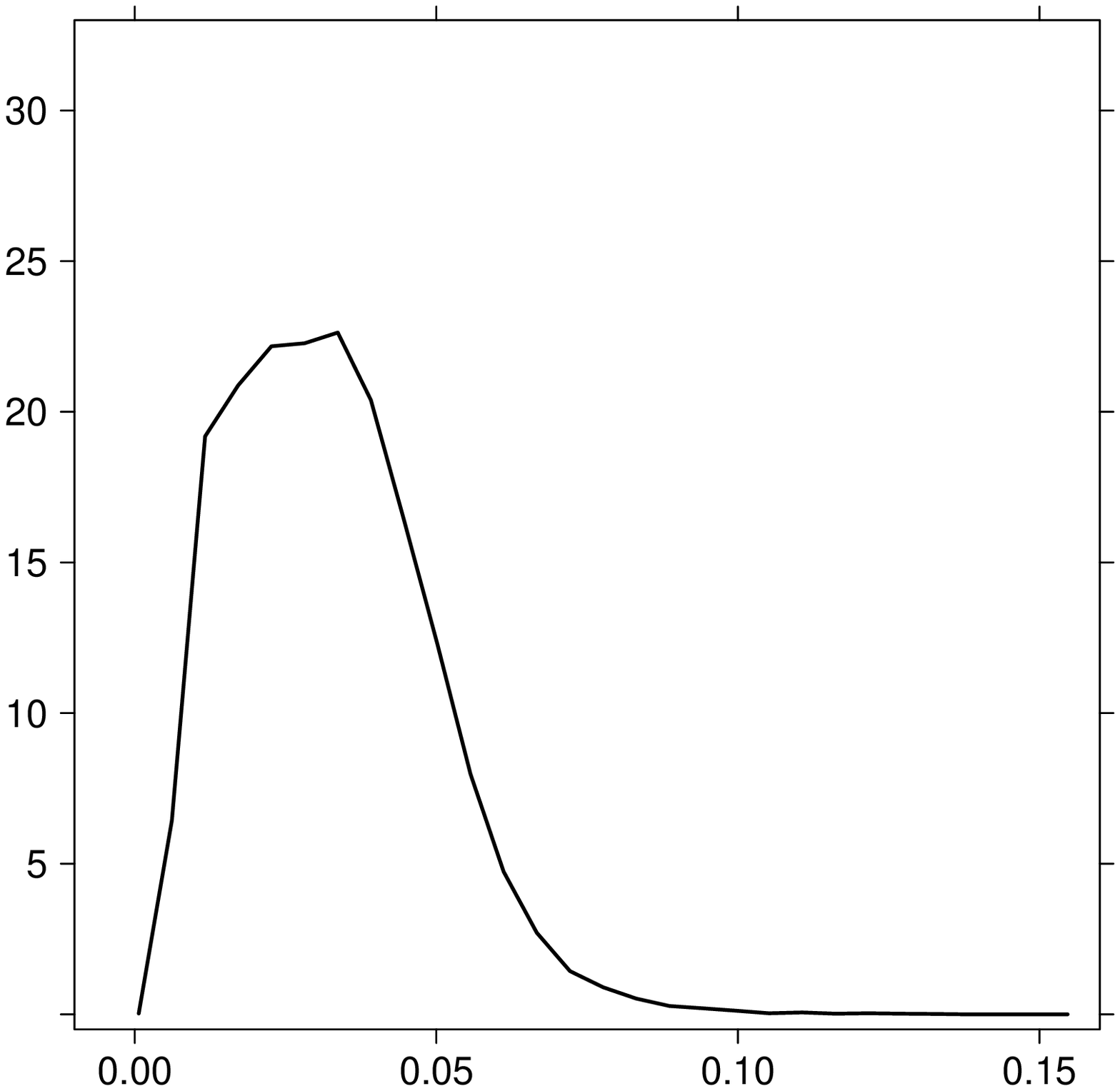}}\hfill\hfill
  \caption{\takep{Posterior credence density of lower bound for
      $\pc^*$.  Data equal to zero have been excluded}\addp{For each
      of the priors of \figref{priors}, the posterior credence density
      of the lower bound for $\pc^*$, conditional on this being
      greater than $0$.}}
  \label{fig:low}
\end{figure}

The lower bound on $\pc^*$ in \eqref{ineqpc*}, $\max\{0,1 - {\Pr(E = 0
  \cd R = 1)}/{\Pr( E=0)}\}$, depends also on $\theta = \Pr( E=1)$,
and its posterior credence distribution could be sensitive to the
prior credence distribution chosen for $\theta$.  We have already
noted that the posterior credence probability that the lower bound is
$0$ is \takep{$0.627$ for Prior~1, and $0.677$ for
  Prior~2.}\addp{$0.65$, independent of the prior for $\theta$.}
\figref{low} displays the posterior densities for the lower bound,
conditional on its being strictly positive, for \takep{these two
  priors}\addp{Prior~1 and Prior~2}; the means are $0.039$ and
$0.025$, and the standard deviations are $0.015$ and $0.016$,
respectively.  We see that the effects of the differences between the
priors are relatively minor.

% \begin{table}[hbtp]
%   \centering
%   \begin{tabular}[c]{lc||c|c}
%     %     \multicolumn{2}{c}{} &
%     &&   \multicolumn{2}{l}{If not $0$:}\\
%     & $\Pr(0)$     & Mean & Standard deviation\\
%     \hline    \hline
%     {Prior~1:} & 0.627 & 0.039 & 0.015\\
%     \hline
%     {Prior~2:} &0.677 & 0.025 & 0.016\\
%     \hline
%   \end{tabular}
%   \caption{Summaries of posterior credence distribution of lower bound on $\pc^*$}
%   \label{tab:lower}
% \end{table}

\subsubsection{Length of interval}
\label{sec:length} 
\begin{figure}[htbp]
  \centering \hfill \subfigure[Prior~1]
  {\includegraphics[scale=.3]{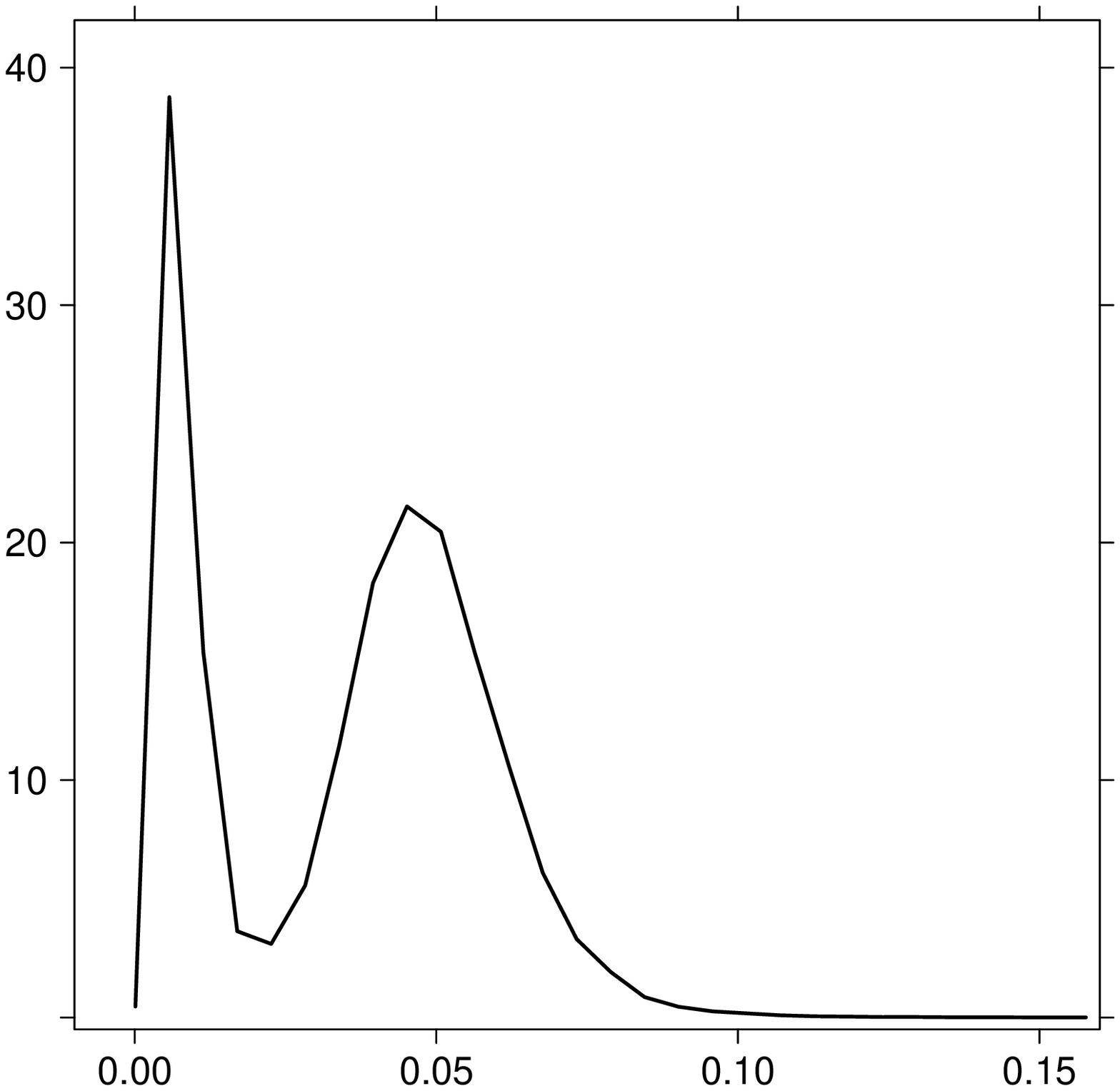}}\hfill \subfigure[Prior~2]
  {\includegraphics[scale=.3]{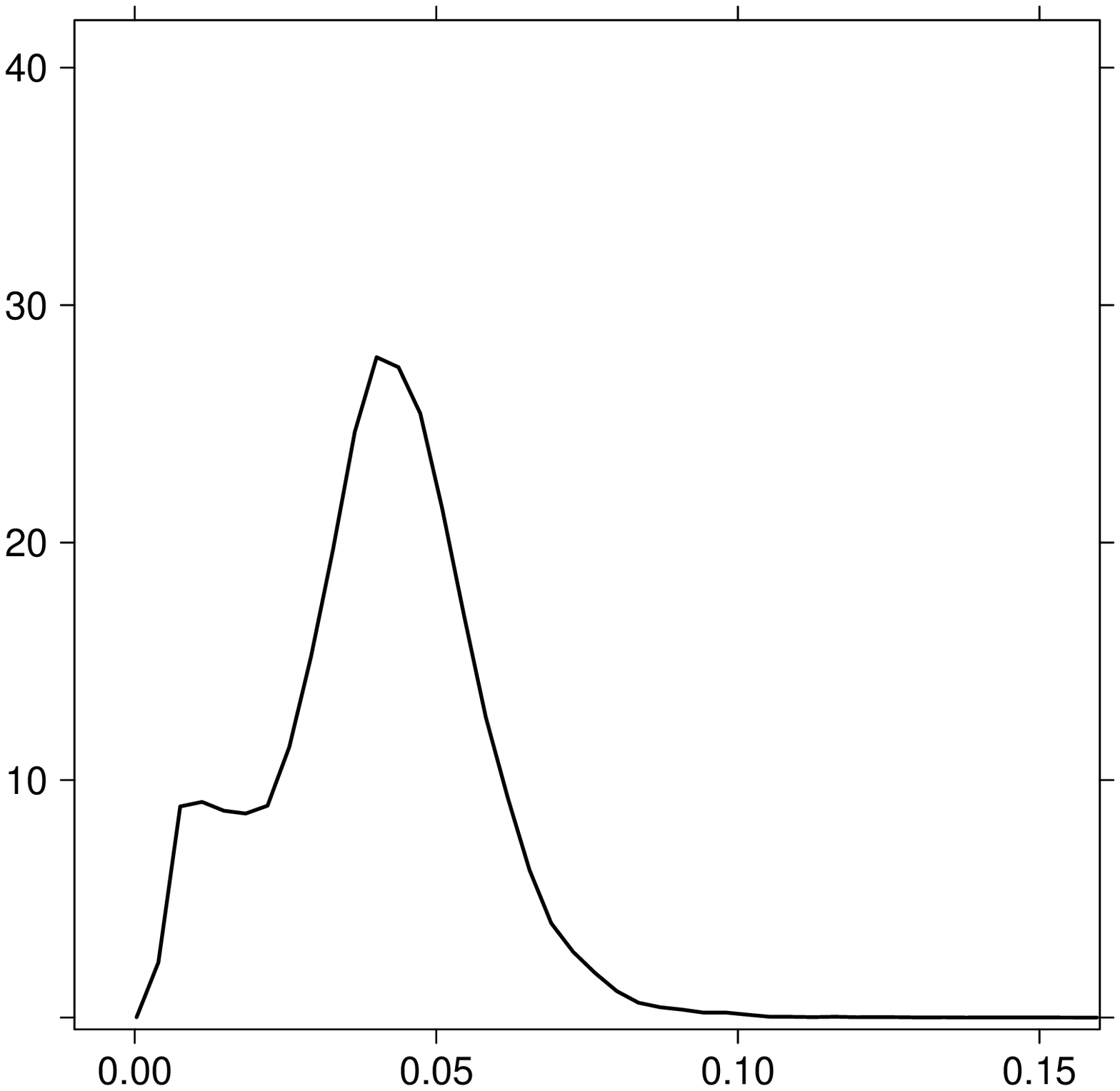}}\hfill\hfill
  \caption{\takep{Posterior credence density of length of interval for
      $\pc^*$}\addp{For each of the priors of \figref{priors}, the
      posterior credence density of the length of \addp{the
        uncertainty} interval for $\pc^*$.}}
  \label{fig:Int}
\end{figure}

Another useful summary of the full inference is the
posterior credence distribution of the length of the interval between
the lower and upper bounds on $\pc^*$, as displayed in \figref{Int}.

\begin{figure}[htbp]
  \centering \hfill \subfigure[Prior~1]
  {\includegraphics[scale=.3]{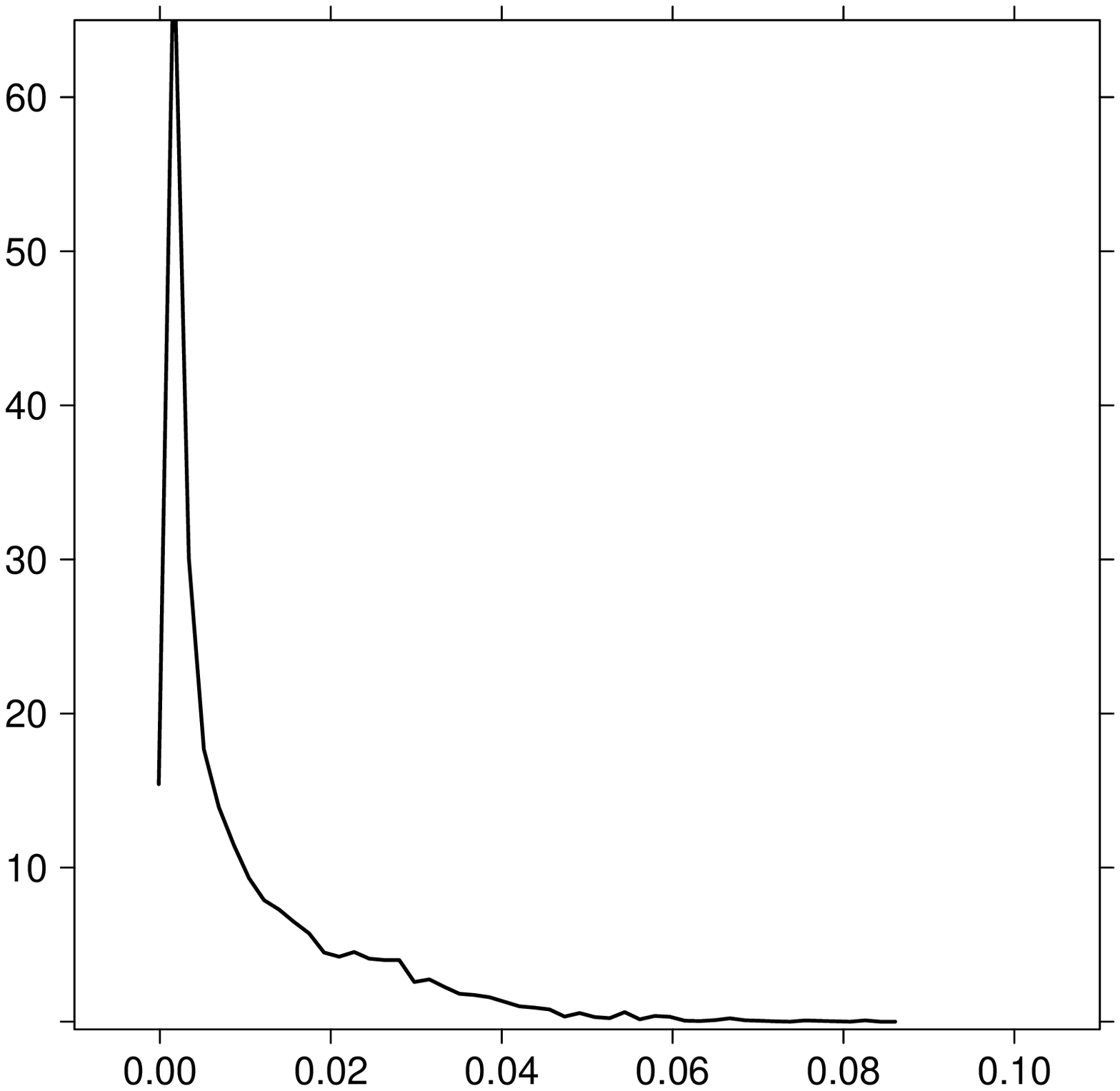}}\hfill \subfigure[Prior~2]
  {\includegraphics[scale=.3]{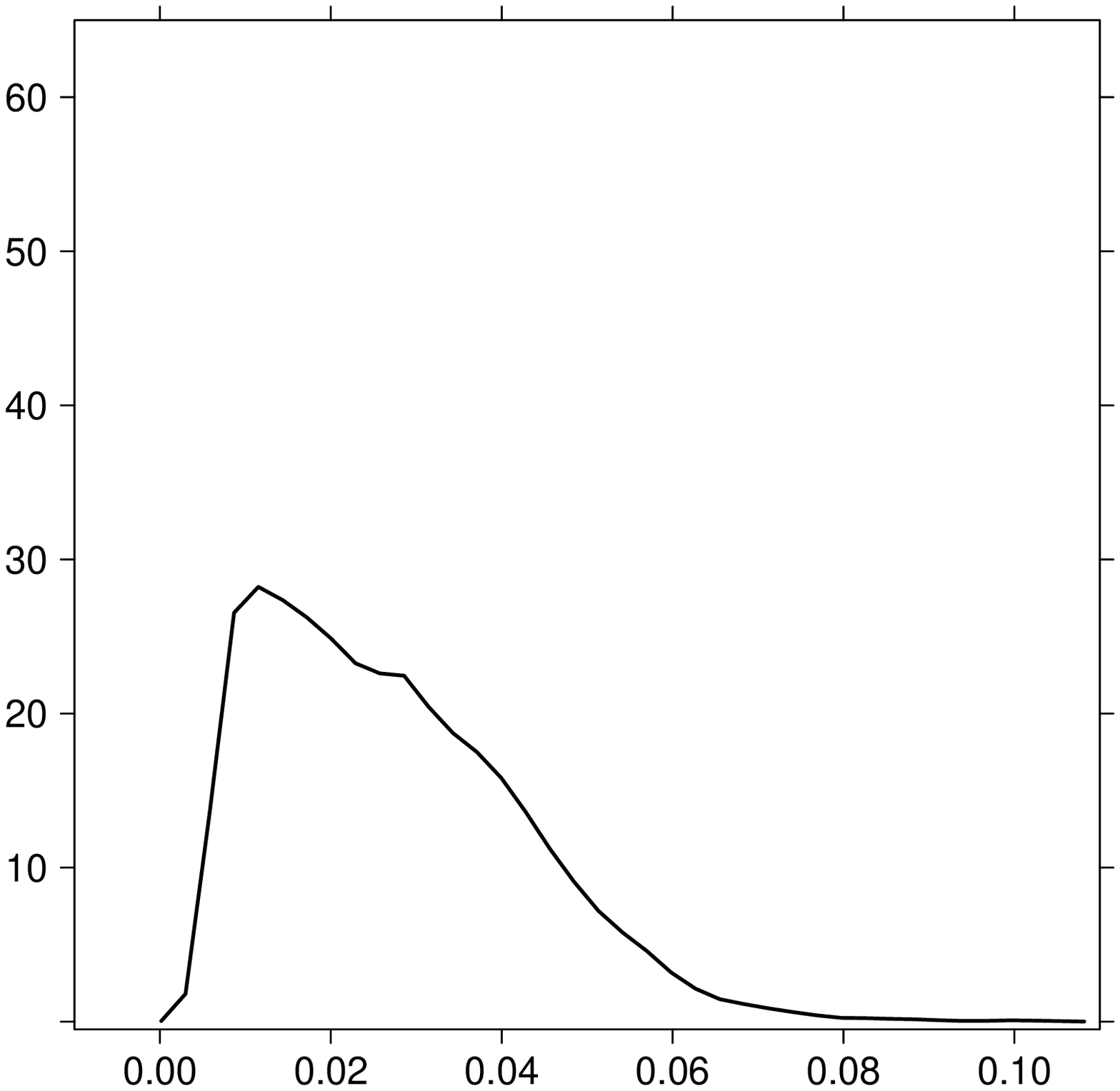}}\hfill\hfill
  \caption{\addp{For each of the priors of \figref{priors}, the
      p}\takep{P}posterior credence density of \addp{the} length of
    \addp{the uncertainty} interval for $\pc^*$, \takep{for data
      with}\addp{conditional on the} lower bound \takep{different from
      zero}\addp{being greater than $0$.}}
  \label{fig:Int1}
\end{figure}

The posterior mean and standard deviation based on Prior~1 are,
respectively, $0.028$ and $0.022$, while for Prior~2 these quantities
are $0.035$ and $0.016.$ We see high sensitivity to the prior
assumptions.  This is particularly apparent when we exclude data with
lower bound equal to zero (see \figref{Int1}).
% Prior~1 has a higher credence for a very short interval
% (corresponding to high precision in the fix on $\pc^*$), at the same
% time as having a higher credence for a long (and so imprecise)
% interval.
For cases with lower bound equal to $0$, the interval length is
identical with the upper bound, as displayed in \figref{UB0}, \takep{,
  where the two priors give very similar results.}\addp{and is
  independent of the prior distribution for $\theta$.} These features
are also \takep{clearly} visible in \figref{lengths}.

% \begin{figure}[htbp]
%   \centering \hfill \subfigure[Prior~1]
%   {\includegraphics[scale=.3]{IP1.eps}}\hfill \subfigure[Prior~2]
%   {\includegraphics[scale=.3]{IP2.eps}}\hfill\hfill
%   \caption{Posterior credence density of length of interval for
%   $\pc^*$. }
%   \label{fig:lengths}
% \end{figure}

\subsubsection{Coverage probability}
\label{sec:cover} 
\begin{figure}[htbp]
  \centering \hfill \subfigure[Prior~1]
  {\includegraphics[scale=.3]{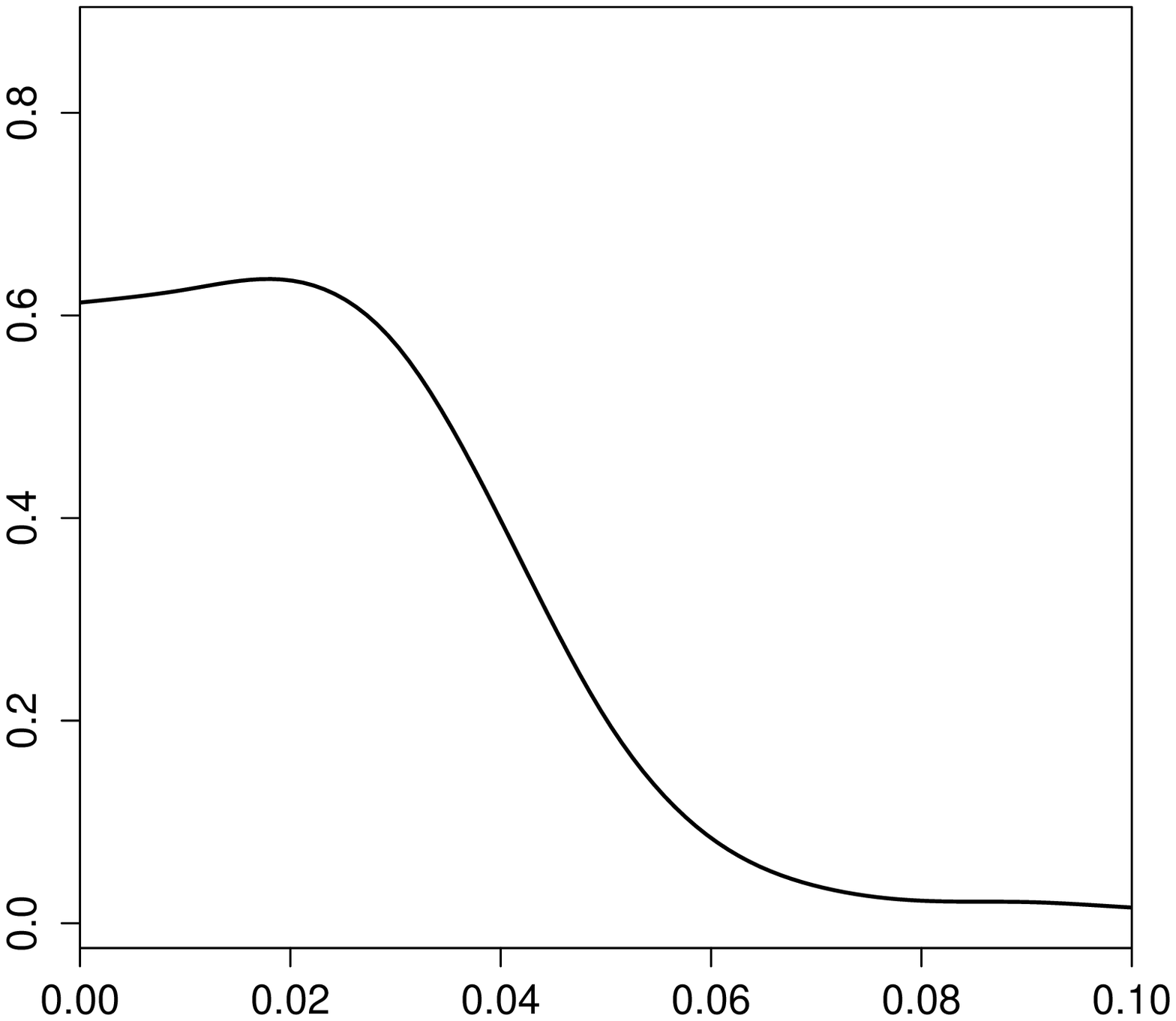}}\hfill \subfigure[Prior~2]
  {\includegraphics[scale=.3]{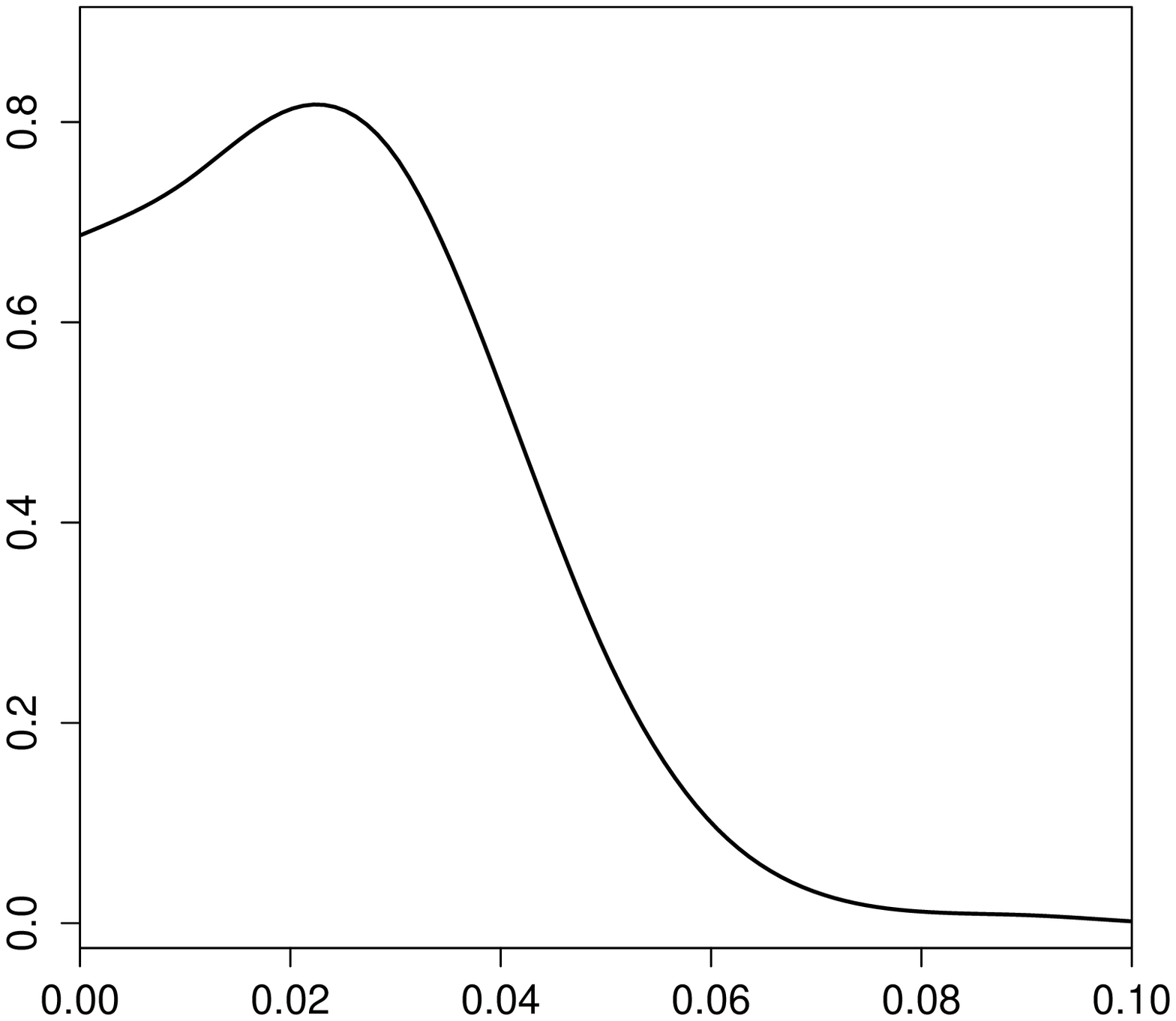}}\hfill\hfill
  \caption{\addp{For each of the priors of \figref{priors}, the
      p}\takep{P}osterior credence probabi\addp{l}ity that the
    \addp{random uncertainty} interval covers a\addp{ny} specific
    value\addp{.}}
  \label{fig:cover}
\end{figure}

Finally, for any probability value $p$, we can
compute the posterior credence that this is included in the random
interval \eqref{ineqpc*}---and thus is at least a candidate as a value
for $\pc^*$.  We graph this coverage measure, as a function of $p$, in
\figref{cover} for both priors.
% Again we note a high degree of sensitivity to prior assumptions.

\addp{
  \subsection{Individual-focused inference}
  \label{sec:indiv}

  The individual-focused inference is much simpler in form: according
  to the analysis in \secref{altview} (and assuming the approximation
  mentioned there is valid), we simply replace the chances featuring
  in the bounds of \eqref{ineqpc*} by their posterior expectations.
  (Recall that we are taking $H$ to be trivial, so it can be omitted
  from the notation).

  The posterior expectation of the upper bound, $\Pr(E = 1 \cd R=1)$,
  is $0.043$, independent of the assumed prior distribution for
  $\theta$.

  As for the lower bound, the posterior expectation of $\Pr(E = 0 \cd
  R = 1)$ is $1-0.043 = 0.957$.  Also, $\Pr( E=0)=1-\theta$, and since
  we have no data relevant to $\theta$ the posterior expectation of
  this quantity is the same as its prior expectation, namely $0.5$ for
  Prior~1, or $0.9$ for Prior~2.  It is clear that there could be high
  sensitivity to the prior distribution assessed for $\theta$.
  However, in this case the lower bound is $0$ for both priors.  Hence
  our individual-focused uncertainty interval for $\pc^*$ is $(0,
  0.043)$ in both cases.

}

\notep{Needs further discussion}

\section{Conclusions}
\label{sec:conc}

We have seen that statistical inference about ``causes of effects'' is
particularly problematic from many points of view, and difficult to
justify even in ideal circumstances.

First, in order merely to formalise the question, we need to carefully
specify, separately, both who is making the inference (in
\secref{statcoe} we called that person``I'') and who (there called
``Ann'') the inference relates to.  Next, we need to be satisfied that
my information $H$ about Ann is {\em sufficient\/}, in the sense of
there being no confounding that could make Ann's treatment choice
informative (for me) about her potential outcome variables.  When all
these conditions are satisfied we can begin to try and learn from
relevant data about the two versions, $\pc$ and $\pc^*$, of the
probability of causation.  For that purpose we should have good
experimental data from which we can get good estimates of the
distribution of the outcome, conditional on exposure $E$ and $H$.  And
even with such ideal estimated probabilities, the resulting inferences
are complex, compounding as they do three different kinds of
uncertainty: interval bounds, for a probability, that are themselves
random.  We have made a start at exploring ways of understanding,
describing and displaying such triple uncertainty (in an example that
admittedly falls far short of the ideal situation), but much remains
to be done.\\

\takep{In the case study of \secref{datanal} we addressed the question
  whether the \alte and \bleed were together caused by \abuse.  But we
  can formulate other CoE questions, such as whether the \alte alone
  was caused by \abuse.  Since we have observed \bleed, this would
  involve replacing the denominator of the lower bound in
  \eqref{ineqpc*} by $1-\Pr(\abuse \cd \bleed)$.  Again, we have no
  data directly relevant to $\Pr(\abuse \cd \bleed)$, and would need
  to assess a prior credence distribution for it.  This might
  reasonably be taken to be higher and tighter than that for
  $\Pr(\abuse)$ alone, but it would again be important to investigate
  sensitivity to a range of reasonable choices.}

\noindent{\bf Acknowledgement.}  We thank Catherine Laurent for
valuable background about Ir\`ene Frachon and benfluorex.

\bibliographystyle{oupvar} \bibliography{ecitsstatcaus,newrefs}
\end{document}